\documentclass[preprint,12pt]{elsarticle}
\usepackage{times}
\usepackage{amsmath,amssymb}
\usepackage{color}
\usepackage{epsfig}
\usepackage{subfigure}
\usepackage{amsfonts}
\usepackage{graphicx}
\usepackage{CJK}
\usepackage{amssymb}
\usepackage{amsthm}
\usepackage{natbib}
\usepackage{amsmath}
\usepackage{setspace}

\allowdisplaybreaks

\newcommand{\pf}{\noindent {\bf Proof. \hspace{2mm}}}

\newcommand{\be}{\begin{equation}}
\newcommand{\ee}{\end{equation}}
\newcommand{\bea}{\begin{eqnarray}}
\newcommand{\eea}{\end{eqnarray}}

\def\pa{\partial}

\allowdisplaybreaks

\begin{document}
 \footskip=0pt
 \footnotesep=2pt
\let\oldsection\section
\renewcommand\section{\setcounter{equation}{0}\oldsection}
\renewcommand\thesection{\arabic{section}}
\renewcommand\theequation{\thesection.\arabic{equation}}
\newtheorem{claim}{\noindent Claim}[section]
\newtheorem{theorem}{\noindent Theorem}[section]
\newtheorem{lemma}{\noindent Lemma}[section]
\newtheorem{proposition}{\noindent Proposition}[section]
\newtheorem{definition}{\noindent Definition}[section]
\newtheorem{remark}{\noindent Remark}[section]
\newtheorem{corollary}{\noindent Corollary}[section]
\newtheorem{example}{\noindent Example}[section]

\title{Regularity of Weak Solutions of Elliptic and Parabolic Equations with Some Critical or Supercritical Potentials}
\author{Zijin Li\footnote{Department of Mathematics and IMS, Nanjing University, Nanjing 210093, China.
 (b081110049@smail.nju.edu.cn)}, Qi S. Zhang\footnote{Department of Mathematics, University of California Riverside, CA 92521, USA. (qizhang@math.ucr.edu)}}

\begin{abstract}
We prove H\"{o}lder continuity of weak solutions of the uniformly elliptic and parabolic  equations
\begin{equation}\label{01}
\partial_{i} ( a_{ij}(x) \partial_{j}u(x)) - \frac{A}{|x|^{2+\beta}} u(x) =0\quad (A>0,\quad\beta\geq 0),
\end{equation}
\begin{equation}\label{02}
\partial_{i} ( a_{ij}(x,t) \partial_{j}u(x,t)) - \frac{A}{|x|^{2+\beta}} u(x,t)-\partial_{t}u(x,t) =0\quad (A>0,\quad\beta\geq 0),
\end{equation}
with critical or supercritical 0-order term coefficients which are beyond  De Giorgi-Nash-Moser's Theory. We also prove, in some special cases, weak solutions are even differentiable.

 Previously P. Baras and J. A. Goldstein \cite{Baras1984} treated the case when $A<0$, $(a_{ij})=I$ and $\beta=0$ for which  they show that there does not exist any regular positive solution or singular positive solutions, depending on the size of $|A|$.  When $A>0$, $\beta=0$ and $(a_{ij})=I$, P. D. Milman and Y. A. Semenov \cite{Milman2003}\cite{Milman2004}
 obtain bounds for the heat kernel.
\end{abstract}


\date{}
\maketitle
\section{Introduction}
In this paper, we consider regularity of weak solutions of divergence form elliptic equations
\begin{equation}\label{maineq1}
\partial_{i} ( a_{ij} \partial_{j}u) - \frac{A}{|x|^{2+\beta}} u=0
\end{equation}
and parabolic equation
\begin{equation}\label{maineq2}
\partial_{i} ( a_{ij} \partial_{j}u) - \frac{A}{|x|^{2+\beta}} u-\partial_{t}u=0
\end{equation}
in the unit ball $B:=B(0,1)\,(or\,\,B\times\mathbb{R}_{+})$ in $\mathbb{R}^{d}$, with $d\geq 3,\,A>0,\,\beta\geq 0$. Here $a_{ij}\in L^{\infty}(B)\,\,(or\,\,L^{\infty}(\mathbb{R}_{+},L^{\infty}(B)))$, and the second order coefficient matrix $\Big(a_{ij}\Big)_{1\leq i,j\leq d}$ satisfies the uniformly elliptic condition:
\begin{equation}
\lambda I\leq \Big(a_{ij}\Big)_{1\leq i,j\leq d}\leq\Lambda I,\text{for some}\,\,0<\lambda\leq\Lambda<\infty.
\end{equation}
Here and below, we use the Einstein summation convention. We say $u\in H^{1}(B)$  is a weak solution of the elliptic equation of \eqref{maineq1}, if $\forall\,\,\psi\in C^{\infty}_{0}(B)$, there holds
\begin{equation}
\int_{B}a_{ij}(x)\partial_{i}\psi(x)\partial_{j}u(x)dx+\int_{B}\frac{A}{|x|^{2+\beta}}u(x)\psi(x)dx=0,
\end{equation}
where $\partial_{i}$ indicates $\partial_{x_i}$ here and below. Similarly, for the parabolic equation \eqref{maineq2}, we say $u\in L^{2}([0,T],H^{1}_{0}(B))$ is a weak solution, if  $\forall\,\,\psi\in C^{\infty}_{0}(B\times[-T,T])$ , there holds
\begin{equation}
\begin{split}
-\int_{0}^{T}\int_{B}u(x,t)\partial_{t}\psi(x,t)dx+\int_{0}^{T}\int_{B}a_{ij}(x,t)\partial_{i}\psi(x,t)\partial_{j}u(x,t)dxdt\\
+\int_{0}^{T}\int_{B}\frac{A}{|x|^{2+\beta}}u(x,t)\psi(x,t)dxdt=\int_{B}\psi(x,0)u(x,0)dx.
\end{split}
\end{equation}



In the middle of the last century, De Giorgi\cite{DeGiorgi1957}, Nash\cite{Nash1958} and Moser\cite{Moser1960}\cite{Moser1961} developed new methods on the studying of elliptic and parabolic equation, which opened  a new area on the study of regularity of weak solutions of elliptic and parabolic equations in divergence form:
\begin{equation}\label{divergenceeeq}
\partial_{i}(a^{ij}\partial_{j}u)+b^{i}\partial_{i}u+cu=f+\partial_{i}f^{i},
\end{equation}
\begin{equation}\label{divergencepeq}
\partial_{i}(a^{ij}\partial_{j}u)+b^{i}\partial_{i}u+cu-\partial_{t}u=f+\partial_{i}f^{i}.
\end{equation}

They proved that, under certain integrable conditions of the coefficients $b^{i}$, $c$ and non-homogeneous term $f$ and $f^{i}$, weak solutions of equation \eqref{divergenceeeq},\eqref{divergencepeq} have $C^{\alpha}$ H\"{o}lder continuity. A key condition for their theory for elliptic equation is that the coefficient of  the 0-order term $c$ must belong to the Lebesgue space $L^{p}$, with $p>\frac{d}{2}$.  Obviously, the 0-order terms in our equations \eqref{maineq1} and \eqref{maineq2} do not satisfy this assumption. Actually the case when $\beta=0$ is the critical borderline case where the theory of De Giorgi, Nash and Moser fails. See for example of Baras and Goldstein \cite{Baras1984} in the case $A<0$.
When $A>0$, even though it is easily seen that weak solutions are locally bounded, it is not clear these
weak solutions have any regularity.

However, such equations are closely related to several physical equations.
For instance the  3-dimensional axially-symmetric incompressible Navier-Stokes equations in fluid dynamics
are
\begin{equation}\label{NS}
\left\{
\begin{array}{l}
\partial_{t}v^{r}+b\cdot\nabla v^{r}-\frac{(v^{\theta})^{2}}{r}+\partial_{r}p=\Big(\Delta-\frac{1}{r^{2}}\Big)v^{r}\\
\partial_{t}v^{\theta}+b\cdot\nabla v^{\theta}+\frac{v^{r}v^{\theta}}{r}=\Big(\Delta-\frac{1}{r^{2}}\Big)v^{\theta}\\
\partial_{t}v^{z}+b\cdot\nabla v^{z}+\partial_{z}p=\Delta v^{z}\\
b=v^{r}r_{r}+v^{z}e_{z},\quad\nabla\cdot b=\partial_{r}v^{r}+\frac{v^{r}}{r}+\partial_{z}v^{z}=0\\
\end{array}
\right.
\end{equation}
where
\begin{equation}
v(x,t)=v^{r}(r,z,t)e_{r}+v^{\theta}(r,z,t)e_{\theta}+v^{z}(r,z,t)e_{z}.
\end{equation}
Observe that the linear parts of the first and second equation of \eqref{NS} are related to equation \eqref{maineq2} with $\beta=0$.

We point out here that the case when $A<0$ was studied in \cite{Baras1984} Baras and Goldstein, who proved that the Cauchy problem of heat equation
\begin{equation}\label{Baras}
\left\{
\begin{array}{l}
\Delta u(x,t)-\frac{A}{|x|^{2}}u(x,t)-\partial_{t}u(x,t)=0,\quad (x,t)\in\mathbb{R}^{d}\times\mathbb{R}_{+}\\
u(x,0)=u_{0}\\
\end{array}
\right.
\end{equation}
have no weak nonnegative solution if $-A>\Big(\frac{d-2}{2}\Big)^{2}$. They also prove if $0<-A\leq\Big(\frac{d-2}{2}\Big)^{2}$, \eqref{Baras} has unbounded positive weak solutions.

The case when $A>0$, $\beta=0$ and the leading operator being the Laplacian was first studied in
Milman and Semenov \cite{Milman2003},
where the authors obtained a sharp upper bound for the fundamental solution of \eqref{maineq2}. Their method is
to use explicit special solutions of the elliptic equation as weights and convert the studying of the
problem to that of a weighted equation via Doob's transform.  Based on this bound, in this special case, one can prove H\"older
continuity of solutions quite easily. Here we observe further that when $A$ is sufficiently large, weak solutions are even differentiable.
  In the variable coefficient case that we
are working on, it is hard or impossible to find an explicit solution of the elliptic equation. So a different method is needed to prove H\"older continuity of weak solutions.

The following are the main results of the paper. The first one pertains the elliptic equation in the case when the leading
operator is the Laplacian but the result is stronger, including differentiability in some situations. This is a little unexpected since
it is well known that singular potential terms usually mess up the derivative bound for solutions.
We also obtain a similar result for the corresponding parabolic cases with $\beta=0$.
The second result deals with both elliptic
and parabolic equations whose leading coefficients are just bounded. We prove H\"{o}lder continuity
of weak solutions.

\begin{theorem}\label{Thm1}
A weak solution $u=u(x)\in H_{0}^{1}(B)$ of $\Delta u-\frac{A}{|x|^{2+\beta}}u=0$ has the following regularity properties. Let $\alpha=\alpha(A)=\frac{-d+2+\sqrt{d^{2}-4d+4+4A}}{2} \in ]n, n+1]$ for a nonnegative integer $n$.


(I) If $\beta=0,$ then $u\in C^{n, (\alpha(A)-n)_-}\big(B_{1/4}\big)$;

(II) If $\beta>0$, then $u\in C^{\infty}\big(B_{1/4}\big)$.
\par\par
\leftline{In addition,}

(III) A weak solution $u=u(x,t)\in L^{2}([0,T],H^{1}_{0}(B))$ of $\Delta u-\frac{A}{|x|^{2}}u-\partial_{t}u=0$ satisfies $\partial_{t}^{m_1}\nabla_{x}^{m_2}u\in C^{(\alpha(A)-n)_-;(\alpha(A)-n)_-/2}\big(B_{1/4}\times[t_{0},T]\big)$, for $2m_{1}+m_{2}=n$.

Here, the H\"{o}lder norms above depend on $d$, $A$ and the $L^{2}$ norm of $u$. $T$, $t_0$ are given
positive constants, and $C_{-}$ defines any number smaller than but close to the constant $C$.
\end{theorem}
\rightline{$\Box$}


\begin{theorem}\label{Var}
If $\beta=0$, the weak solution $u$ of the elliptic equation \eqref{maineq1} is H\"{o}lder continuous, i.e.
\begin{equation}
\|u\|_{C^{\alpha}(B_{1/4})}\leq C(\lambda,\Lambda,d,A)\|u\|_{L^{2}(B)}
\end{equation}
with $\alpha=\alpha(\lambda,\Lambda,d,A)>0$. Moreover any weak solution $u=u(x,t)$ of the parabolic equation \eqref{maineq2} is H\"{o}lder continuous when $t$ is away from $0$, i.e.
\begin{equation}
\|u\|_{C^{\alpha;\frac{\alpha}{2}}(B_{1/4}\times[t_{0},T] )}\leq C(\lambda,\Lambda,d,A,)t_{0}^{-(d/2+1+\alpha)}\|u\|_{L^{2}}
\end{equation}
with $\alpha=\alpha(\lambda,\Lambda,d,A)>0$, $t_0>0$.
\end{theorem}
\rightline{$\Box$}
This theorem provides an interior estimate which deteriorates near initial time. However, this is necessary since no H\"older regularity assumption on the initial datum is made.

 We also mention that the Harnack inequality could not hold for  solutions of these equations, because one can find a class of non-negative solutions which do not satisfy it. See section \ref{Lap} e.g..  Moreover, these special solutions are instrumental in studying the regularity of the solutions of
\begin{equation}\label{ceq}
\Delta u-\frac{A}{|x|^{2+\beta}}u=0,\quad(\beta\geq 0),
\end{equation}
the elliptic case where the leading operator is the Laplacian.
 It helps to prove an $\alpha-$order decay estimate of the weak solution $u=u(x)$ when $\beta=0$ at $0\in\mathbb{R}^{d}$, namely:
\begin{equation}
|u(x)|\leq C|x|^{\alpha}
\end{equation}
for $\alpha\in(0,1)$ and $C\in\mathbb{R}_{+}$ is a constant. When $\beta>0$, the decay of $u(x)$ at $0\in\mathbb{R}^{d}$ turns exponential.

For the variable 2nd order coefficients case \eqref{maineq1}, \eqref{maineq2}, the situation is more complicated. Roughly speaking, we could not find a good enough special solution as the Laplacian case \eqref{ceq}. However, if $\beta=0$, we find a weighted mean value inequality, which is motivated by \cite{Zhang2000} and \cite{Wong2003}. The weight, decaying at certain rate near the origin, plays the same role as the special solution in the Laplacian case \eqref{ceq}, giving similar $\alpha-$decay estimate.



The rest of the paper is organized as follows.  In section 2, we give the proof of Theorem 1.1. In section 3, we state and prove the aforementioned weighted mean value inequality
for general parabolic equations \eqref{maineq2}. In section 4, we give the proof of the variable coefficient case for weak solutions of \eqref{maineq1} with critical 0-order term coefficient. Finally in section 5, we extend our conclusion in section 4 to the parabolic case. Some elementary but useful works, giving the proof of the existence of weak solutions, local boundedness of weak solution, maximum principle, and an introduction of the modified Bessel's equation, could be found in the Appendix.


\section{Laplacian Case}\label{Lap}

In this section, we will prove Theorem \ref{Thm1}. First we need a simple
lemma on certain special solutions of \eqref{ceq}, which will serve as a benchmark for
comparison with other solutions.  As mentioned in the introduction, the case when $\beta=0$, H\"older
continuity of solutions
can also be proven by the bound in \cite{Milman2003}. Here we give a direct proof based on the maximal
principle.

\begin{lemma}\label{Prop1}
(i) If $\beta=0$, then
\begin{equation}
u(x)=|x|^{\alpha}, \quad \alpha=\alpha(A)=\frac{-d+2+\sqrt{d^{2}-4d+4+4A}}{2}
\end{equation}
is a weak solution of \eqref{ceq}.\\
(ii) If $\beta>0$, then
\begin{equation}
u(x)=|x|^{-\frac{d}{2}+1}\mathcal{K}_{(d-2)/\beta}\big(\frac{2}{\beta}\sqrt{A}|x|^{-\frac{\beta}{2}}\big)
\end{equation}
is a weak solution of \eqref{ceq}, where $\mathcal{K}_{(d-2)/\beta}$ is the \emph{modified Bessel's function of second kind} mentioned above.
\end{lemma}
\pf
Since we are looking for radially symmetric solution of \eqref{ceq} here, we can just solve the corresponding
ODE.  Define $r=|x|$, we find the solution $u=u(r)$. Thus, \eqref{ceq} turns to\\
\begin{equation}\label{1}
u''+\frac{d-1}{r}u'-\frac{A}{r^{2+\beta}}u=0
\end{equation}
If $\beta=0$, this is an Euler type ODE. Set $u=r^{\alpha}$, and take this into the equation \eqref{1}, we have
\begin{equation}
\alpha^{2}+(d-2)\alpha-A=0
\end{equation}
Solve this equation with a positive number, we have\begin{equation}\alpha=\alpha(A)=\frac{-d+2+\sqrt{d^{2}-4d+4+4A}}{2}.
\end{equation}

If $\beta>0$, suppose $\mathcal{B}_{\lambda}(r)$ satisfies the modified Bessel's equation
\begin{equation}
r^{2}\mathcal{B}_{\lambda}''(r)+r\mathcal{B}_{\lambda}'(r)-(r^{2}+\lambda^{2})\mathcal{B}_{\lambda}(r)=0.
\end{equation}
Change $r$ into $\nu\cdot r^{\mu}$, where $\nu\neq 0$ and $\mu\neq 0$ are real numbers to be determined later, we have $g(r):=\mathcal{B}_{\lambda}(\nu\cdot r^{\mu})$ satisfies the following equation by direct calculation:
\begin{equation}\label{2}
g''+\frac{1}{r}g'-\frac{\mu^{2}\nu^{2}}{r^{2-2\mu}}g-\frac{\mu^{2}\lambda^{2}}{r^{2}}g=0.
\end{equation}
Observe \eqref{1}, we choose $\mu=-\frac{\beta}{2}$ and $\nu=\frac{2}{\beta}\sqrt{A}$, thus \eqref{2} becomes
\begin{equation}\label{3}
g''+\frac{1}{r}g'-\frac{A}{r^{2+\beta}}g-(\frac{\beta\lambda}{2r})^{2}g=0.
\end{equation}
Now, to eliminate the last term and modify the coefficient of the second term on the left of \eqref{3}, we set $h(r):=r^\theta\cdot g(r)$, where $\theta$ is a real number to be determined later. By direct calculation, we have $h$ satisfies
\begin{equation}\label{4}
h''+\frac{1-2\theta}{r}h'-\frac{A}{r^{2+\beta}}h+\big[\theta^{2}-\big(\frac{\beta\lambda}{2}\big)^{2}\big]r^{-2}h=0.
\end{equation}
Compare \eqref{4} to \eqref{1}, we have
\begin{equation}
\left\{
\begin{array}{l}
\theta^{2}-\big(\frac{\beta\lambda}{2}\big)^{2}=0\\
 \\
-2\theta+1=d-1\\
\end{array}
\right.
\end{equation}
Thus we have $\theta=-\frac{d}{2}+1$, and $\lambda=\frac{d-2}{\beta}$, which we have $h(r)=r^{-\frac{d}{2}+1}\mathcal{B}_{\frac{2-d}{\beta}}(\frac{2}{\beta}\sqrt{A}r^{-\frac{\beta}{2}})$ solves \eqref{1}. Since we are looking for local bounded solution, we choose
\begin{equation}
\mathcal{B}_{\frac{d-2}{\beta}}(r)=\mathcal{K}_{\frac{d-2}{\beta}}(r),
 \end{equation}
 \emph{the modified Bessel's function of second kind}, which is exponentially growing at $0\in\mathbb{R}^{d}$.

\rightline{$\Box$}
As for (ii) in Lemma \ref{Prop1} , we have

\begin{lemma}
The function $|x|^{-\frac{d}{2}+1}\mathcal{K}_{(d-2)/\beta}\big(\frac{2}{\beta}\sqrt{A}|x|^{-\frac{\beta}{2}}\big)$ is smooth in $B$, and decays exponentially to 0 at $x=0$.
\end{lemma}
This is a direct corollary of the property of modified Bessel's function, see \cite{Abramowitz1972} for more details.\\
\rightline{$\Box$}

Now we start the proof of Theorem \ref{Thm1},  {\bf Case (I)}.

We denote by $J_\beta=J_{\beta}(x)$ the special solutions of \eqref{ceq} mentioned above, namely
\begin{equation}
J_{\beta}(x)=\left\{
\begin{array}{l}
|x|^{-\frac{d}{2}+1}\mathcal{K}_{(d-2)/\beta}\big(\frac{2}{\beta}\sqrt{A}|x|^{-\frac{\beta}{2}}\big),\quad \beta>0\\
 \\
|x|^{\frac{-d+2+\sqrt{d^{2}-4d+4+4A}}{2}},\quad \beta=0.\\
\end{array}
\right.
\end{equation}
Then, on $B_{1/2}$, there exists a constant $C=C\big(d,\|u\|_{L^{2}(B)}\big)$, such that the functions $v_1=u(x) + CJ_{\beta}(x)$ and  $v_2=u(x) - C J_{\beta}(x)$ satisfy
\begin{equation}\label{aueq}
\left\{
\begin{array}{l}
\Delta v_i(x)-\frac{1}{|x|^{2+\beta}}v_i(x)=0,\quad in\,\,B_{1/2}, \quad i=1, 2,\\
 \\
v_1(x)\geq 0,\quad\,\,v_2(x)\leq 0,\quad on\,\,\partial B_{1/2}.\\
\end{array}
\right.
\end{equation}
By the maximum principle in Lemma \ref{Wmp} in the Appendix e.g., we have
\begin{equation}\label{ubd}
-C\cdot J_{\beta}(x)\leq u(x)\leq C\cdot J_{\beta}(x).
\end{equation}

According to the Green's Representation formula(c.f.\cite{GT1998}), we have, for $x\in B_{1/2}$
\begin{equation}\label{30}
\begin{split}
u(x)&=\int_{B_{1/2}}\Gamma(x-y)\frac{1}{|y|^{2+\beta}}u(y)dy+H(x)\\
\end{split}
\end{equation}
where $\Gamma(x-y)=\frac{|x-y|^{2-d}}{d(2-d)\omega_{d}}$ is the fundamental solution of the Laplace equation ($\omega_d$ is the volume of d-dimensional unit ball), and $H=H(x)$ is harmonic in $B_{1/2}$. 
Since the second term $H(x)$ of \eqref{30} is regular enough in $B_{1/4}$,  we only need to consider the regularity of
\begin{equation}
w(x)=\int_{B_{1/2}}\Gamma(x-y)\frac{1}{|y|^{2+\beta}}u(y)dy.
\end{equation}
We divide the rest of the proof into several cases, firstly:
\subsection{\textbf{Case (I) $\beta=0,\,\,\alpha(A)\leq 1$}}
By \eqref{ubd},$\forall x_{1},x_{2}\in B_{1/4}$, define $z=\frac{1}{2}(x_{1}+x_{2}),\,\,\delta=|x_{1}-x_{2}|$.
\begin{equation}
\begin{split}
|w(x_{1})-w(x_{2})|&\leq C\int_{B_{1/2}}|\Gamma(x_{1}-y)-\Gamma(x_{2}-y)|\frac{1}{|y|^{2-\alpha(A)}}dy\\
&\leq C\int_{B_{1/2}\bigcap B(z,\delta)}|\Gamma(x_{1}-y)|\frac{1}{|y|^{2-\alpha(A)}}dy\\
&+C\int_{B_{1/2}\bigcap B(z,\delta)}|\Gamma(x_{2}-y)|\frac{1}{|y|^{2-\alpha(A)}}dy\\
&+C\int_{B_{1/2}- B(z,\delta)}|\Gamma(x_{1}-y)-\Gamma(x_{2}-y)|\frac{1}{|y|^{2-\alpha(A)}}dy\\
&=C(I_{1}+I_{2}+I_{3})\\
\end{split}
\end{equation}
As to $I_{1}$, for $\frac{d}{2}<p_{1}<\frac{d}{2-\alpha(A)}$, by H\"{o}lder inequality
\begin{equation}\label{I1}
\begin{split}
I_{1}&\leq C\Big(\int_{B(x_{1},\frac{3\delta}{2})}|x_{1}-y|^{\frac{(2-d)p_{1}}{p_{1}-1}}dy\Big)^{1-1/p_{1}}\cdot\Big(\int_{B_{1/2}}|y|^{(\alpha(A)-2)p_{1}}dy\Big)^{1/p_{1}}\\
&\leq C \delta^{\frac{2p_{1}-d}{p_{1}}}\cdot\Big[\frac{1}{(\alpha(A)-2)p_{1}+d}\Big]^{\frac{1}{p_{1}}}\cdot\Big[\frac{2p_{1}-d}{p_{1}-1}\Big]^{-1+1/p_{1}}\\
&\leq C\cdot\delta^{\alpha(A)_{-}},\quad \text{ by choosing}\,\,p_{1}\to\big(\frac{d}{2-\alpha(A)}\big)_{-}.
\end{split}
\end{equation}
Here and below, we use $C_{-}$ to denote an arbitrary number close but smaller than $C$. Similarly, we have $I_{2}$ satisfies the same estimate. As for $I_{3}$, by mean value inequality, there exists an $\hat{x}$ lies between $x_{1}$ and $x_{2}$, and for $1<p_{2}<\frac{d}{2-\alpha(A)}$:
\begin{equation}
\begin{split}
I_{3}&\leq C\delta\int_{B_{1/2}- B(z,\delta)}|\nabla\Gamma(\hat{x}-y)|\frac{1}{|y|^{2-\alpha(A)}}dy\\
& \leq C\delta\Big(\int_{B(0,\frac{\delta}{2})^{c}}|y|^{\frac{(1-d)p_{2}}{p_{2}-1}}dy\Big)^{1-1/p_{2}}\cdot\Big(\int_{B_{1/2}}|y|^{(\alpha(A)-2)p_{2}}dy\Big)^{1/p_{2}}\\
&\leq C \delta^{\frac{2p_{2}-d}{p_{2}}}\cdot\Big[\frac{1}{(\alpha(A)-2)p_{2}+d}\Big]^{\frac{1}{p_{2}}}\\
&\leq C\cdot \delta^{\alpha(A)_{-}},\quad \text{by choosing}\,\,p_{2}\to\big(\frac{d}{2-\alpha(A)}\big)_{-}.
\end{split}
\end{equation}
This means $w=w(x)$ is $\alpha(A)_{-}$ H\"{o}lder continuous in $B_{1/4}$. This and
\eqref{30} imply that
\[
\Vert  u \Vert_{C^{\alpha(A)_{-}}\left(B_{1/4}\right)} \le C(d, \Vert u \Vert_{L^2(B)}).
\] This proves Case (I) of the theorem when $n=0$.

We point out that this estimate is almost optimal since the special solution $J_{0}(x)=|x|^{\alpha(A)}$ has only $\alpha(A)$ H\"{o}lder continuity in $B_{1/2}$. \rightline{$\Box$}
\begin{remark}
Let us pay attention to a special case when $A=1$, $\beta=0$, $d=3$. According to the results above, we know that the weak solution of
\begin{equation}
\Delta u-\frac{1}{|x|^{2}}u=0
\end{equation}
in $B_{1/2}$ is H\"older continuous with exponent $\big(\frac{\sqrt{5}-1}{2}\big)_{-}\approx 0.618$, which is the golden ratio.
\end{remark}

Next we prove  Case (I) of the theorem when $n=1$.
In this case, we first claim $w\in C^{1}\big(B_{1/4}\big)$
 and
\begin{equation}\label{1stD}
\partial_{i}w(x)=\int_{B_{1/2}}\partial_{i}\Gamma(x-y)\frac{u(y)}{|y|^{2}}dy,\quad i=1,2,...,d.
\end{equation}
Here goes the proof.

Since $\nabla^{n}\Gamma$ satisfies the following estimate
\begin{equation}
|\nabla^{n}\Gamma(x-y)|\leq C|x-y|^{2-d-n},\quad n=1,2,...
\end{equation}
where $C=C(d,n)$,  the following function
\begin{equation}
\xi(x)=\int_{B_{1/2}}\partial_{i}\Gamma(x-y)\frac{u(y)}{|y|^{2}}dy
\end{equation}
is well defined by the H\"{o}lder inequality. The reason is:
\begin{equation}
\begin{split}
\partial_{i}\Gamma(x-y)\in L^{p}(B_{1/2})&,\quad 1\leq p<\frac{d}{d-1};\\
\frac{|u(y)|}{|y|^{2}}\leq C|y|^{\alpha(A)-2}\in L^{q}(B_{1/2})&,\quad 1\leq q<\frac{d}{2-\alpha(A)}.
\end{split}
\end{equation}
Thus
\begin{equation}
1-\frac{\alpha(A)-1}{d}<\frac{1}{p}+\frac{1}{q}\leq 2
\end{equation}
and we get
\begin{equation}
\partial_{i}\Gamma(x-y)\frac{|u(y)|}{|y|^{2}}\in L^{(1+\frac{\alpha(A)-1}{d-\alpha(A)+1})_{-}}(B_{1/2}).
\end{equation}
Therefore $\xi$ is well defined. By usual approximation argument,  one can prove easily that $\partial_{i}w(x)=\xi(x)$. This proves the claim.

Similarly as in \textbf{Case (I)}, $n=0$, $\forall x_{1},x_{2}\in B_{1/4}$, define $z=\frac{1}{2}(x_{1}+x_{2}),\,\,\delta=|x_{1}-x_{2}|$. We have the  inequality for the gradient of $w$:
\begin{equation}
\begin{split}
|\partial_{i}w(x_{1})-\partial_{i}w(x_{2})|&\leq C\int_{B_{1/2}}|\partial_{i}\Gamma(x_{1}-y)-\partial_{i}\Gamma(x_{2}-y)|\frac{1}{|y|^{2-\alpha(A)}}dy\\
&\leq C\int_{B_{1/2}\bigcap B(z,\delta)}|\partial_{i}\Gamma(x_{1}-y)|\frac{1}{|y|^{2-\alpha(A)}}dy\\
&+C\int_{B_{1/2}\bigcap B(z,\delta)}|\partial_{i}\Gamma(x_{2}-y)|\frac{1}{|y|^{2-\alpha(A)}}dy\\
&+C\int_{B_{1/2}- B(z,\delta)}|\partial_{i}\Gamma(x_{1}-y)-\partial_{i}\Gamma(x_{2}-y)|\frac{1}{|y|^{2-\alpha(A)}}dy\\
&\equiv C(I_{1}+I_{2}+I_{3}).\\
\end{split}
\end{equation}
As to $I_{1}$, for $d<p_{1}<\frac{d}{2-\alpha(A)}$, by H\"{o}lder inequality
\begin{equation}\label{I1}
\begin{split}
I_{1}&\leq C\Big(\int_{B(x_{1},\frac{3\delta}{2})}|x_{1}-y|^{\frac{(1-d)p_{1}}{p_{1}-1}}dy\Big)^{1-1/p_{1}}\cdot\Big(\int_{B_{1/2}}|y|^{(\alpha(A)-2)p_{1}}dy\Big)^{1/p_{1}}\\
&\leq C \delta^{\frac{p_{1}-d}{p_{1}}}\cdot\Big[\frac{1}{(\alpha(A)-2)p_{1}+d}\Big]^{\frac{1}{p_{1}}}\cdot\Big[\frac{2p_{1}-d}{p_{1}-1}\Big]^{-1+1/p_{1}}\\
&\leq C\cdot\delta^{(\alpha(A)-1)_{-}},\quad \text{ by choosing}\,\,p_{1}\to\big(\frac{d}{2-\alpha(A)}\big)_{-}.
\end{split}
\end{equation}
Likewise, we have that $I_{2}$ satisfies the same estimate. As for $I_{3}$, by mean value inequality, there exists an $\hat{x}$ lies between $x_{1}$ and $x_{2}$, and for $1<p_{2}<\frac{d}{2-\alpha(A)}$:
\begin{equation}
\begin{split}
I_{3}&\leq C\delta\int_{B_{1/2}- B(z,\delta)}|\nabla^{2}\Gamma(\hat{x}-y)|\frac{1}{|y|^{2-\alpha(A)}}dy\\
& \leq C\delta\Big(\int_{B(0,\frac{\delta}{2})^{c}}|y|^{-\frac{dp_{2}}{p_{2}-1}}dy\Big)^{1-1/p_{2}}\cdot\Big(\int_{B_{1/2}}|y|^{(\alpha(A)-2)p_{2}}dy\Big)^{1/p_{2}}\\
&\leq C \delta^{1-\frac{d}{p_{2}}}\cdot\Big[\frac{1}{(\alpha(A)-2)p_{2}+d}\Big]^{\frac{1}{p_{2}}}\\
&\leq C\cdot \delta^{(\alpha(A)-1)_{-}},\quad \text{by choosing}\,\,p_{2}\to\big(\frac{d}{2-\alpha(A)}\big)_{-}.
\end{split}
\end{equation}
This means $\nabla w$ is $(\alpha(A)-1)_{-}$ H\"{o}lder continuous in $B_{1/4}$ and $w\in C^{1,\,\,(\alpha(A)-1){-}}\big(B_{1/4}\big)$.  Moreover
\begin{equation}
\Vert \nabla u \Vert_{C^{(\alpha(A)-1)_{-}}\big(B_{1/4}\big)} \le C(d, \Vert u \Vert_{L^2(B)}).
\end{equation}
This shows Case (I) of the theorem when $n=1$ holds.

Now we prove {\bf Case (I)} of the theorem when $n>1$.
First by induction, it is easy to see the following lemma, we list here without proof.
\begin{lemma}
If $u\in C^{m,\gamma}\left(B\right)$, $\gamma\in(0,1)$ and $|u(x)|\leq C|x|^{n+\gamma}$ with $n\geq m$, we have
\begin{equation}\label{230}
|\nabla^{k}u(x)|\leq  C(\|u\|_{C^{m,\gamma}\left(B\right)}) \cdot|x|^{n-k+\gamma},\forall k\in[1,m]\text{ and } k \text{ is an integer,}
\end{equation}
for all $x\in B$.
\end{lemma}
\qed

Moreover,
\begin{equation}\label{2.32}
\begin{split}
&\hskip 1cm \partial_{x_j} \int_{B_{1/2}}\partial_{x_i}\Gamma(x-y)\nabla^{n-2}\left(\frac{u(y)}{|y|^{2}}\right)dy\\
&=-P.V. \int_{B_{1/2}} \partial_{x_j} \partial_{y_i}\Gamma(x-y)\nabla^{n-2}\left(\frac{u(y)}{|y|^{2}}\right)dy\\
&=\int_{B_{1/2}} \partial_{x_j} \Gamma(x-y)\partial_{y_i}\nabla^{n-2}\left(\frac{u(y)}{|y|^{2}}\right)dy-\int_{\partial B_{1/2}} \partial_{x_j} \Gamma(x-y)\nabla^{n-2}\left(\frac{u(y)}{|y|^{2}}\right) n_i dS_{y}.
\end{split}
\end{equation}
Note the last term above is smooth since $x_1$, $x_2\in B_{1/4}$.  By \eqref{30}, in order to prove H\"older continuity of $\nabla^n u$, we only need to prove that of
\begin{equation}
w_n(x):=\int_{B_{1/2}} \partial_{x_j} \Gamma(x-y)\partial_{y_i}\nabla^{n-2}\left(\frac{u(y)}{|y|^{2}}\right)dy
\end{equation}
for each $n$. Just as the situation $n=0$ and $1$, since $u$ and its derivatives satisfy the decay property \eqref{230}, we get the following regularity result by induction.
\begin{equation}
\Vert \nabla^{n} w \Vert_{C^{(\alpha(A)-n)_{-}}\big(B_{1/4}\big)} \le C(d, \Vert u \Vert_{L^2(B)}).
\end{equation}Consequently
\begin{equation}
\Vert \nabla^{n} u \Vert_{C^{(\alpha(A)-n)_{-}}\big(B_{1/4}\big)} \le C(d, \Vert u \Vert_{L^2(B)}).
\end{equation}
Thus we finished the proof of Case I of  Theorem 1.

We also mention here that we could use the classical Schauder estimate to get the regularity of higher derivatives, by treating the term $A u(x)/|x|^2$ as the inhomogeneous term and using the decay property
of $u$ near $0$.

For Case (II), i.e., $\beta>0$, one can use \eqref{ubd} and the exponential decay property of
$J_\beta$ to conclude that $u$ decays exponentially at $0$. Then following an analogous argument
as in Case (I), we know that $u \in C^\infty(B_{1/4})$.

Before commencing the case of variable second order coefficients, we prove part (III) of the Theorem \ref{Thm1}, namely the
 corresponding result for the parabolic equations.  The result is based on  a bound of the fundamental solution $\Gamma_{0}=\Gamma_{0}(x,t;y,s)$ of  the parabolic equation with leading term being Laplace
\begin{equation}\label{CPB}
\Delta u-\frac{A}{|x|^{2}}u-\pa_{t}u=0,
\end{equation}
which was proved in \cite{Milman2004}. It states that
\begin{equation}\label{2.37}
\Gamma_{0}(x,t;y,0)\leq \frac{C}{t^{d/2}}\left(1+\frac{\sqrt{t}}{|x|}\right)^{-\alpha}\left(1+\frac{\sqrt{t}}{|y|}\right)^{-\alpha}e^{-c\frac{|x-y|^{2}}{t}}.
\end{equation}
Where $C$ and $c$ are positive constants, $\alpha=\alpha(A)=\frac{-d+2+\sqrt{d^{2}-4d+4+4A}}{2}$, just as before. See also \cite{Alexander2006} for more details. We claim here that this bound of fundamental solution leads to a higher regularity of weak solutions of \eqref{CPB} in a neighbourhood of $x=0$ when $t$ is away from zero, 
by a similar method of Lemma 2.2 in \cite{Wong2003}, also Lemma \ref{Mean} in the next section, we have a mean value inequality of $u$:
\begin{equation}\label{CMeanValue}
u^{2}(x,t)\leq \frac{C}{|Q_{r}(x,t)|}\left(\frac{|x|}{r}\right)^{2\alpha}\int_{Q_{2r}(x,t)}u^{2}(y,s)dyds,
\end{equation}
where $C>0$ and $Q_{r}(x,t)=B(x,r)\times [t-r^2,t]\subset\mathbb{R}^{d}\times\mathbb{R}_{+}$. This inequality gives an "$\alpha-$order" decay of $u$ when $t$ is away from 0, namely:
\begin{equation}\label{Pbound}
|u(x,t)|\leq C(d,\|u\|_{L^{2}}) t^{-(d/2+1+\alpha)}|x|^{\alpha}.
\end{equation}
Now we are going to consider the regularity of $u=u(x,t)$ with $t>t_{0}>0$. By Duhamel's Principle (see also section 5 for more details), we only need to consider the regularity of
\begin{equation}
w(x,t)\equiv\int_{t_0/2}^{t}\int_{B_{1/2}}\frac{e^{-\frac{|x-y|^{2}}{4(t-s)}}}{(t-s)^{d/2}}\frac{u(y,s)}{|y|^{2}}dyds.
\end{equation}
Since $\alpha\in]0,1]$ will be considered in section 5, we only consider $\alpha\in]n,n+1]$ for an integer $n\geq 1$. Firstly, if $n=1$, similarly as \eqref{1stD}, we have
\begin{equation}\label{1stDP}
\begin{split}
\partial_{x_i}w(x,t)=&\int_{t_0/2}^{t}\int_{B_{1/2}}\partial_{x_{i}}\left(\frac{e^{-\frac{|x-y|^{2}}{4(t-s)}}}{(t-s)^{d/2}}\right)\frac{u(y,s)}{|y|^{2}}dyds\\
=&-\frac{1}{2}\int_{t_{0}/2}^{t}\int_{B_{1/2}}\frac{e^{-\frac{|x-y|^{2}}{4(t-s)}}(x_{i}-y_{i})}{(t-s)^{d/2+1}}\frac{u(y,s)}{|y|^{2}}dyds.
\end{split}
\end{equation}
Using the same method in section 5 (from \eqref{MEstimate} to \eqref{(5.22)}), we have the H\"older continuity of $\nabla_{x} w$, hence we have:
\begin{equation}
\|\nabla_{x} u\|_{C^{(\alpha(A)-1)_{-};(\alpha(A)-1)_{-}/2}\left( B_{1/4}\times[t_{0},T]\right)}\leq C\left(d,t_0,\|u\|_{L^{2}}\right).
\end{equation}
This shows the result holds when $n=1$. When $n>1$, by a similar induction method as in the elliptic case (see \eqref{2.32} for more detials), we have the H\"older continuity of higher order spacial derivatives, namely:
\begin{equation}
\|\nabla_{x}^{n} u\|_{C^{(\alpha(A)-n)_{-};(\alpha(A)-n)_{-}/2}\left(B_{1/4}\times[t_{0},T] \right)}\leq C\left(d,t_0,\|u\|_{L^{2}}\right).
\end{equation}
Finally, we will estimate time derivatives of $u$. By equation \eqref{CPB}, $\partial_{t}u=\Delta u-\frac{A}{|x|^{2}}u$, we can eliminate all the time derivatives, for $2m+|L|=n$
\begin{equation}\label{2.44}
\partial_{t}^{m}\partial_{x}^{L}u=\sum_{|J|+|K|+2l=n}C_{J,K,l,d}\frac{x^{J}}{|x|^{2(l+|J|)}}\left(\partial_{x}^{K}u\right)(x,t).
\end{equation}
Where $J$, $K$, $L$ above are multi-indexes, e.g. $L=(L_{1}, L_{2},...,L_{d})$, $|L|=L_{1}+L_{2}+...+L_{d}$
\begin{equation}
\partial_{x}^{L}=\partial_{x_1}^{L_1}\partial_{x_2}^{L_2}...\partial_{x_d}^{L_d};\quad x^{L}=x_{1}^{L_1}x_{2}^{L_2}...x_{d}^{L_d}.
\end{equation}
And $C_{J,K,l,d}$ is a constant. Then we will give a H\"older estimate of each term on the right hand side of \eqref{2.44}.
We choose $x_{1}, x_{2}\in B_{1/4}$ and $T>t_{1}>t_{2}>t_{0}>0$. We suppose $|x_{1}|\geq|x_{2}|$ without loss of generality. Then
\begin{equation}
\begin{split}
&\left|\frac{x_{1}^{J}}{|x_{1}|^{2(l+|J|)}}\partial_{x}^{K}u(x_{1},t_{1})-\frac{x_{2}^{J}}{|x_{2}|^{2(l+|J|)}}\partial_{x}^{K}u(x_{2},t_{2})\right|\\
\leq&\left|\frac{x_{1}^{J}}{|x_{1}|^{2(l+|J|)}}\partial_{x}^{K}u(x_{1},t_{1})-\frac{x_{1}^{J}}{|x_{1}|^{2(l+|J|)}}\partial_{x}^{K}u(x_{1},t_{2})\right|\\
&+\left|\frac{x_{1}^{J}}{|x_{1}|^{2(l+|J|)}}\partial_{x}^{K}u(x_{1},t_{2})-\frac{x_{1}^{J}}{|x_{1}|^{2(l+|J|)}}\partial_{x}^{K}u(x_{2},t_{2})\right|\\
&+\left|\frac{x_{1}^{J}}{|x_{1}|^{2(l+|J|)}}\partial_{x}^{K}u(x_{2},t_{2})-\frac{x_{2}^{J}}{|x_{2}|^{2(l+|J|)}}\partial_{x}^{K}u(x_{2},t_{2})\right|\\
:=&I+II+III
\end{split}
\end{equation}
We only consider cases when $|K|<n$, otherwise $\frac{x^{J}}{|x|^{2(l+|J|)}}\partial_{x}^{K}u=\partial_{x}^{K}u=\partial_{x}^{n}u$ which was considered before. Actually, by the structure of equation \eqref{CPB}, we have $|K|\leq n-2$ and $\partial_{t}\partial^{K}_{x}u(x_{1},\cdot)$ is well defined according to the equation \eqref{CPB}. Then, if $t_{1}-t_{2}\leq|x_{1}|^{2}$,
\begin{equation}
\begin{split}
I&\leq \frac{C}{|x_{1}|^{2l+|J|}}\sup_{\tau\in[t_{0},T]}|\partial_{t}\partial_{x}^{K}u(x_{1},\tau)|\cdot(t_{1}-t_{2})\\
&\leq \frac{C}{|x_{1}|^{2l+|J|}}|x_{1}|^{n+2-\alpha(A)}\left(\sup_{\tau\in[t_{0},T]}|\Delta\partial_{x}^{K}u(x_{1},\tau)|+\sup_{\tau\in[t_{0},T]}\left|\frac{\partial_{x}^{K}u(x_{1},\tau)}{|x_{1}|^{2}}\right|\right)\cdot(t_{1}-t_{2})^{(\alpha(A)-n)/2}\\
&\leq\frac{C(t_{0})}{|x_{1}|^{2l+|J|}}\cdot|x_{1}|^{n+2-\alpha(A)}\cdot|x_{1}|^{\alpha(A)-(|K|+2)}(t_{1}-t_{2})^{(\alpha(A)-n)/2}\\
&=C(t_{0})(t_{1}-t_{2})^{(\alpha(A)-n)/2}.
\end{split}
\end{equation}
If $t_{1}-t_{2}>|x_{1}|^{2}$,
\begin{equation}
\begin{split}
I&\leq \frac{C}{|x_{1}|^{2l+|J|}}\cdot\frac{(t_{1}-t_{2})^{(\alpha(A)-n)/2}}{|x_{1}|^{\alpha(A)-n}}\left(|\partial_{x}^{K}u(x_{1},t_{1})|+|\partial_{x}^{K}u(x_{1},t_{2})|\right)\\
&\leq \frac{C(t_{0})}{|x_{1}|^{2l+|J|}}\cdot\frac{(t_{1}-t_{2})^{(\alpha(A)-n)/2}}{|x_{1}|^{\alpha(A)-n}}\cdot|x_{1}|^{\alpha(A)-|K|}\\
&=C(t_{0})(t_{1}-t_{2})^{(\alpha(A)-n)/2}.
\end{split}
\end{equation}
Thus we finished the estimate of I. As for II and III, we have, for a $\hat{x}$ lies between $x_1$ and $x_2$,
\begin{equation}
\begin{split}
II&\leq\frac{C}{|x_{1}|^{2l+|J|}}|\nabla_{x}\partial_{x}^{K}u(\hat{x},t_{2})|\cdot|x_{1}-x_{2}|\leq\frac{C|\hat{x}|^{\alpha(A)-1-|K|}}{|x_{1}|^{2l+|J|}}|x_{1}-x_{2}|\\
&\leq\frac{C}{|x_{1}|^{n-\alpha(A)+1}}(|x_{1}|+|x_{2}|)^{n-\alpha(A)+1}\cdot|x_{1}-x_{2}|^{\alpha(A)-n}\leq C2^{n-\alpha(A)+1}\cdot|x_{1}-x_{2}|^{\alpha(A)-n};
\end{split}
\end{equation}
\begin{equation}
\begin{split}
III&\leq|\partial_{x}^{K}u(x_{2},t_{2})|\cdot\left|\frac{x_{1}^{J}}{|x_{1}|^{2(l+|J|)}}-\frac{x_{2}^{J}}{|x_{2}|^{2(l+|J|)}}\right|\leq C\cdot|x_{2}|^{\alpha(A)-|K|}\left|\frac{x_{1}^{J}}{|x_{1}|^{2(l+|J|)}}-\frac{x_{2}^{J}}{|x_{2}|^{2(l+|J|)}}\right|\\
&\leq C\cdot|x_{1}-x_{2}|^{\alpha(A)-n}.
\end{split}
\end{equation}

This proves that, for $2m+|L|=n$:
\begin{equation}
\|\partial_{t}^{m}\partial_{x}^{L}u\|_{C^{(\alpha(A)-n)_{-};(\alpha(A)-n)_{-}/2}\left( B_{1/4}\times[t_{0},T]\right)}\leq C\left(d,t_0,\|u\|_{L^{2}}\right).
\end{equation}
\qed

We can also draw the same conclusion by using the fact that $\partial_{t}^{n}u$ is also solution of \eqref{CPB} in some sense. But certain approximation procedure is needed since we do not know a priori
any $L^p$ bounds for the time derivatives.

\section{A Mean Value Inequality}

Now we start treating equations with variable coefficients.

In this section, we state and prove a mean value inequality for solutions of equations
\eqref{02}, which has extra decay comparing with the standard mean value inequality. This will be used in the following sections. The proof uses an iteration process involving the potential and a boosting process by the Feynman-Kac formula. To start with, we need a crude mean value inequality, which is similar to that in \cite{Wong2003}
\begin{lemma}\label{Crude}
Let $u$ be a weak solution to the equation
\begin{equation}\label{93}
\partial_{x_{i}} ( a_{ij} \partial_{x_{j}}u) - V u-\partial_{t}u =0
\end{equation}
in $Q_{2r}(x,t)$, Here $a_{ij}\in L^{\infty}$, $\lambda I\leq \big(a_{ij}\big)_{1 \leq i,j\leq d}\leq\Lambda I$ for some $0<\lambda\leq\Lambda\leq\infty$. If $a\geq a_{0}:=\frac{192(d+4)^{2}(1+\Lambda)}{\lambda}$ and
\begin{equation}
V=\frac{a}{1+|x|^{2}},
\end{equation}
then $\exists\,\,C>0$, depending on $a$, $\lambda$, $\Lambda$ but independent of $r$, such that:
\begin{equation}
u^{2}(x,t)\leq C\frac{[\max\{r/(1+|x|),1\}]^{-2}}{|Q_{r}(x,t)|}\int_{Q_{r}(x,t)}u^{2}(y,s)dyds.
\end{equation}
\end{lemma}
\pf
We pick a Lipschitz cut-off function $\phi$ such that $\phi(y,s)=1$ if $(y,s)\in Q_{r}(x,t)$, $\phi(y,s)=0$ if $(y,s)\in Q_{\tau r}^{c}(x,t)$, and $|\nabla \phi|\leq 1/((\tau-1)r)$,a.e. $|\partial_{t}\phi|\leq 1/((\tau-1)r)^{2}$,a.e. Using $\phi^{2}u$ as a test function, after routine calculation,
\begin{equation}
\lambda\int|\nabla(\phi u)|^{2}dyds+\int V(y)u^{2}\phi^{2}dyds\leq(1+\Lambda)\int u^{2}\big[|\nabla\phi|^{2}+|\partial_{s}\phi|\big]dyds.
\end{equation}
Therefore:
\begin{equation}
\int_{Q_{\tau r}(x,t)}V(y)u^{2}\phi^{2}dyds\leq\frac{1+\Lambda}{\lambda((\tau-1)r)^{2}}\int_{Q_{\tau r}(x,t)}u^{2}dyds.
\end{equation}
When $y\in B(x,r)$, we have $|y|^{2}\leq 2(|x-y|^{2}+|x|^{2})$ and hence, without loss of generality, we assume $|x|\geq 1$ and $r\geq 1$. Otherwise \eqref{Crude} is the standard mean value inequality.
\begin{equation}
V(y)\geq\frac{a}{3(|x-y|^{2}+|x|^{2})}\geq\frac{a}{3(r^2+|x|^2)}.
\end{equation}
It follows that
\begin{equation}
\int_{Q_{r}(x,t)}u^{2}dyds\leq\frac{3(1+\Lambda)(r^{2}+|x|^{2})}{a\lambda(\tau-1)^{2}r^{2}}\int_{Q_{\tau r}(x,t)}u^{2}dyds.
\end{equation}
For each $r>1$ we take $\tau>1$ s.t.
\begin{equation}
\frac{3(1+\Lambda)(r^{2}+|x|^{2})}{a\lambda(\tau-1)^{2}r^{2}}=\frac{1}{2}.
\end{equation}
This implies
\begin{equation}
\tau r=r+\big[6a^{-1}\lambda^{-1}(1+\Lambda)(r^{2}+|x|^{2})\big]^{1/2}.
\end{equation}
Under such choice of $\tau$, we have
\begin{equation}\label{2.4}
\int_{Q_{r}(x,t)}u^{2}dyds\leq\frac{1}{2}\int_{Q_{\tau r}(x,t)}u^{2}dyds.
\end{equation}
We shall iterate the above inequality according to the formula:
\begin{equation}
\tau_{k+1}=\tau_{k}r_{k}:=r_{k}+\big[6a^{-1}\lambda^{-1}(1+\Lambda)(r_{k}^{2}+|x|^{2})\big]^{1/2}
\end{equation}
with $r_{0}=|x|$. Writing $\mu:=\big(12a^{-1}\lambda^{-1}(1+\Lambda)\big)^{1/2}$, we claim that:
\begin{equation}\label{ini}
r_{k}\leq (1+\mu)^{2k}(1+|x|).
\end{equation}
Obvious \eqref{ini} hold s for $k=1$. Suppose it holds for $k$, then
\begin{equation}
\begin{split}
r_{k+1}&\leq r_{k}+\mu(r_{k}+|x|)=(1+\mu)r_{k}+\mu|x|\\
&\leq (1+\mu)^{2k+1}(1+|x|)+(1+\mu)^{2k+1}\mu(1+|x|)\\
&\leq(1+\mu)^{2k+2}(1+|x|).
\end{split}
\end{equation}
This implies that to reach $r$ from $|x|$ one needs at least
\begin{equation}
k=\frac{\ln(r/(1+|x|))}{2\ln(1+\mu)}
\end{equation}
number of iterations (round up to an integer). Iterating \eqref{2.4} $k$ times we have
\begin{equation}
\begin{split}
u^{2}(x,t)&\leq\frac{C}{|Q_{r_{0}(x,t)}|}\int_{Q_{r_{0}}(x,t)}u^{2}(y,s)dyds\leq\frac{C2^{-k}}{|Q_{r_{0}(x,t)}|}\int_{Q_{r}(x,t)}u^{2}(y,s)dyds\\
&\leq C e^{-\frac{\ln[r/(1+|x|)]}{\ln(1+\mu)}\frac{\ln 2}{2}}\frac{(r/r_{0})^{d+2}}{|Q_{r}(x,t)|}\int_{Q_{r}(x,t)}u^{2}(y,s)dyds.
\end{split}
\end{equation}
Simplifying the above, we reach
\begin{equation}
u^{2}(x,t)\leq C\Big(\frac{r}{1+|x|}\Big)^{-\ln2/[2\ln(1+\mu)]}\Big(\frac{r}{1+|x|}\Big)^{d+2}\frac{1}{|Q_{r}(x,t)|}\int_{Q_{r}(x,t)}u^{2}(y,s)dyds.
\end{equation}
Recall that $\mu=\big(12a^{-1}\lambda^{-1}(1+\Lambda)\big)^{1/2}$. When $a\geq\frac{192(d+4)^{2}(1+\Lambda)}{\lambda}$, we have
\begin{equation}
\mu\leq\frac{1}{4(d+4)}.
\end{equation}
Hence
\begin{equation}
\frac{\ln 2}{2\ln(1+\mu)}\geq\frac{\ln2}{2\mu}=(d+4)2\ln2\geq d+4.
\end{equation}
This shows
\begin{equation}
u^{2}(x,t)\leq C\Big(\frac{r}{1+|x|}\Big)^{-2}\frac{1}{|Q_{r}(x,t)|}\int_{Q_{r}(x,t)}u^{2}(y,s)dyds.
\end{equation}
This proves the lemma.

\rightline{$\Box$}

Based on this crude mean value inequality, using the Feymann-Kac product formula similarly as in Proposition 2.1 of \cite{Wong2003}, denoting $t-s$ by l, we obtain the global bounds of $\Gamma_{1}$, the fundamental solution of \eqref{93}.
\begin{equation}\label{94}
\Gamma_{1}(x,t;y,s)=\Gamma_{1}(x,l;y,0)\leq c_{1}\frac{w(x,l)}{|B(x,\sqrt{l})|}e^{-c_{2}|x-y|^{2}/l},
\end{equation}
for any $y\in Q_{3r/2}(x,t)$. Here $w(x,l)=[\max\{\frac{\sqrt{l}}{1+|x|},1\}]^{-\alpha}$ with $\alpha=\alpha(\lambda,\Lambda, a,d)>0$.
Since the proof is the same, we omit the details here. Using this bound, we will give a refined mean value formula. Notice that there is no restriction on the size of the positive number $A$ comparing with Lemma \ref{Crude}.

\begin{lemma}\label{Mean}
Let $u$ be a weak solution to the parabolic equation
\begin{equation}\label{93}
\partial_{x_{i}} ( a_{ij} \partial_{x_{j}}u) - \frac{A}{1+|x|^{2}} u-\partial_{t}u =0
\end{equation}
with $a_{ij}(x,t)\in L^{\infty}$ satisfying the elliptic condition $\lambda I<\big(a_{ij}(x,t)\big)<\Lambda I$ with $0<\lambda<\Lambda<\infty$ in the parabolic cube $Q_{2r}(x,t)=B(x,2r)\times [t-4r^{2},t]$,$\forall r>0$. Then there exists $C>0$, $\alpha>0$, depending only on $\lambda$, $\Lambda$, $A$, $d$, such that
\begin{equation}
u^{2}(x,t)\leq\frac{C[max\{r/(1+|x|),1\}]^{-2\alpha}}{|Q_{r}(x,t)|}\int_{Q_{r}(x,t)}u^{2}(y,s)dyds
\end{equation}
\end{lemma}

\pf
Select a cut-off function $\eta\in C^{\infty}_{0}(Q_{3r/2}(x,t))$ such that $\eta=1$ in $Q_{r}(x,t)$, $0\leq \eta\leq 1$ and $|\nabla \eta|\leq \frac{C}{r}$, $|\nabla^{2}\eta|\leq\frac{C}{r^{2}}$, $|\eta_{t}|\leq\frac{C}{r^{2}}$. We have $\eta u$ satisfies
\begin{equation}
\left\{
\begin{array}{l}
\partial_{j}\big(a_{ij}\partial_{i}(\eta u)\big)-\frac{A}{1+|x|^{2}}\eta u-\partial_{t}(\eta u)=u\partial_{j}\big(a_{ij}\partial_{i}\eta\big)+2a_{ij}\partial_{j}\eta\partial_{i}u-u\partial_{t}\eta:=f;\\
\eta u(y,s)=0, (y,s)\in\partial B(x,3r/2)\times[t-(3r/2)^{2},t];\\
\eta u(y,t-(3r/2)^{2})=0\\
\end{array}
\right.
\end{equation}
Let $\Gamma_{1}(x,t;y,s)$ be the fundamental solution of $\partial_{j}(a_{ij}\partial_{i}u)-\frac{A}{1+|x|^{2}}u-\partial_{t}u=0$. Then
\begin{equation}\label{Green}
\begin{split}
u(x,t)&=-\int_{Q_{3r/2}(x,t)}\Gamma_{1}(x,t;y,s)\big(a_{ij}\partial_{j}\eta\partial_{i}u-u\partial_{t}\eta\big)(y,s)dyds\\
&+\int_{Q_{3r/2}(x,t)}\partial_{j}\Gamma_{1}(x,t;y,s)\big(a_{ij}u\partial_{i}\eta\big)(y,s)dyds\\
&=I+II.
\end{split}
\end{equation}
For I, we have:
\begin{equation}\label{99}
I^{2}\leq C\int_{Q_{3r/2}(x,t)-Q_{r}(x,t)}\Gamma_{1}^{2}(x,t;y,s)dyds\int_{Q_{3r/2}(x,t)-Q_{r}(x,t)}\Big(\frac{1}{r^{4}}u^{2}+\frac{1}{r^{2}}|\nabla u|^{2}\Big)dyds.
\end{equation}
Since $\partial_{j}(a_{ij}\partial_{i}u)-\partial_{t}u=\frac{A}{1+|x|^{2}}u$, it is well known that
\begin{equation}\label{100}
\int_{Q_{3r/2}(x,t)-Q_{r}(x,t)}|\nabla u|^{2}dyds\leq \frac{C}{r^{2}}\int_{Q_{2r}(x,t)-Q{r/2}(x,t)}u^{2}dyds.
\end{equation}
Combine \eqref{99} and \eqref{100} we have
\begin{equation}\label{101}
I^{2}\leq\frac{C}{r^{4}}\int_{Q_{3r/2}(x,t)-Q_{r}(x,t)}\Gamma_{1}^{2}(x,t;y,s)dyds\int_{Q_{2r}(x,t)}u^{2}dyds.
\end{equation}
As for II:
\begin{equation}
II^{2}\lesssim\frac{1}{r^{2}}\int_{Q_{3r/2}(x,t)-Q_{r}(x,t)}|\nabla \Gamma_{1}(x,t;y,s)|^{2}dyds\cdot\int_{Q_{3r/2}(x,t)-Q_{r}(x,t)}u^{2}dyds.
\end{equation}
Similarly as \eqref{100}, since $\Gamma_{1}(x,t;y,s)$ is a well-defined weak solution of equation \eqref{93} in $Q_{3r/2}(x,t)-Q_{r}(x,t)$, we have:
\begin{equation}
II^{2}\lesssim\frac{1}{r^{4}}\int_{Q_{2r}(x,t)-Q_{r/2}(x,t)}\Gamma_{1}(x,t;y,s)^{2}dyds\cdot\int_{Q_{3r/2}(x,t)-Q_{r}(x,t)}u^{2}dyds.
\end{equation}
By the estimate of I and II above, we arrive:
\begin{equation}
u(x,t)^{2}\lesssim\frac{1}{r^{4}}\int_{Q_{2r}(x,t)-Q_{r/2}(x,t)}\Gamma_{1}(x,t;y,s)^{2}dyds\cdot\int_{Q_{3r/2}(x,t)-Q_{r}(x,t)}u^{2}dyds.
\end{equation}
By \eqref{94}, denoting $t-s$ by l, we obtain
\begin{equation}
\Gamma_{1}(x,t;y,s)=\Gamma_{1}(x,l;y,0)\leq c_{1}\frac{w(x,l)}{|B(x,\sqrt{l})|}e^{-c_{2}|x-y|^{2}/l}.
\end{equation}
Here $w(x,l)=[\max\{\frac{\sqrt{l}}{1+|x|},1\}]^{-\alpha}$ with $\alpha=\alpha(\lambda,\Lambda, A,d)>0$.

When $(y,s)\in Q_{2r}(x,t)-Q_{r/2}(x,t)$ and $t-(r/2)^{2}\leq s\leq t$, we have $0\leq ;\leq(r/2)^{2}\leq|x-y|^{2}$. Hence
\begin{equation}
\Gamma_{1}(x,t;y,s)\leq c_{1}\frac{\max\{\frac{\sqrt{l}}{1+|x|},1\}^{-\alpha}}{l^{d/2}}e^{-(c_{2}r^{2})/(4l)}\leq c\frac{\max\{\frac{r}{1+|x|},1\}^{-\alpha}}{r^{d}}.
\end{equation}
On the other hand, when $(y,s)\in Q_{2r}(x,t)-Q_{r/2}(x,t)$ and $t-(2r)^{2}\leq s\leq t-(r/2)^{2}$,
we have $l=t-s\geq (r/2)^{2}$. Therefore
\begin{equation}
\Gamma_{1}(x,t;y,s)\leq c_{1}\frac{w(x,l)}{l^{d/2}}e^{-c_{2}|x-y|^{2}/l}\leq c\frac{\max\{\frac{r}{1+|x|},1\}^{-\alpha}}{r^{d}}
\end{equation}
Therefore, $\forall(y,s)\in Q_{2r}(x,t)-Q_{r/2}(x,t)$, we have
\begin{equation}
\Gamma_{1}(x,t;y,s)\leq c\frac{\max\{\frac{r}{1+|x|},1\}^{-\alpha}}{r^{d}}
\end{equation}
Substituting this into \eqref{101} to get
\begin{equation}
u(x,t)^{2}\leq\max\{\frac{r}{1+|x|},1\}^{-2\alpha}\frac{C}{r^{2+d}}\int_{Q_{2r}(x,t)}u^{2}dyds
\end{equation}
\rightline{$\Box$}
\section{Variable Second Order Coefficients Case}
In this section we prove Theorem \ref{Var}, the elliptic case. Namely
\begin{equation}
\partial_{i} ( a_{ij}(x) \partial_{j}u(x)) - \frac{A}{|x|^{2}} u(x) =0
\end{equation}
where $a_{ij}\in L^{\infty}(B)$ and satisfies uniformly elliptic condition. Although this elliptic case is a special case of the parabolic one in the next section, we present a proof since it is more transparent.

Obviously, this situation is different from the Laplacian case in the section before, since we could not find a special solution as in section 2. But fortunately, for the critical case ($\beta=0$), we can still prove that weak solutions of \eqref{maineq1} vanish at $0\in \mathbb{R}^{d}$ in the order of $|x|^{\alpha}$ for some $\alpha>0$ by the mean value inequality in Lemma \ref{Mean}. First we need the following two lemmas:
\begin{lemma}\label{MeanValue}
Let $u\in H^{1}(B)$ be a weak solution of the equation
\begin{equation}
\partial_{i}(a_{ij}(x)\partial_{j}u(x))-\frac{A}{|x|^{2}}u(x)=0
\end{equation}
$x_{0}\in B_{1/2}\subset\mathbb{R}^{d}$. Here $a_{ij}\in L^{\infty}(B)$, $\lambda I\leq \big(a_{ij}\big)_{1 \leq i,j\leq d}\leq\Lambda I$ for some $0<\lambda\leq\Lambda\leq\infty$. 
Then there exists an $\alpha>0$, $C>0$ depending only on $\lambda$, $\Lambda$, $A$, $d$ and the $L^{2}(B)$ norm of $u$, such that
\begin{equation}
|u(x_{0})|\leq C|x_{0}|^{\alpha}
\end{equation}
\end{lemma}
\pf
Since $1/|x|^{2}$ is bounded except in $B(0,\delta)$, $\delta>0$ fixed, by standard theory, $u$ is H\"{o}lder continuous on $B(0,\delta)^{c}$. Note that this H\"{o}lder exponent may not be uniform when $\delta\to 0$. Pick a positive integer $k$, $r\geq|x_{0}|\geq 0$ and consider the equation
\begin{equation}\label{Apeq}
\left\{
\begin{array}{l}
\partial_{i}\big(a_{ij}(x)\partial_{j}u_{k}(x)\big)-\frac{A}{|x|^{2}+k^{-2}}u_{k}(x)=0\\
u_{k}\big|_{\partial B(x_{0},2r)}=u\big|_{\partial B(x_{0},2r)}\\
\end{array}
\right.
\end{equation}
 By boundedness of the potentials $A/(|x|^{2}+k^{-2})$, the problem above has a unique solution which is H\"{o}lder continuous on $B(x_{0},2r)-B(0,\delta)$, $\forall\delta>0$ fixed, with H\"{o}lder exponent and H\"{o}lder norms that are uniform with respect to $k$ for each fixed $\delta$. By Arzela-Ascoli theorem, there exists a subsequence, still denoted by $\{u_{k}\}$ for simplicity, such that
\begin{equation}\label{limit}
u_{k}(x)\to u(x)\quad \text{uniformly in } B(x_{0},2r)-B(0,\delta)
\end{equation}
Next we do a scaling $x=y/k$, define $\tilde{u}_{k}(y)=u_{k}(\frac{y}{k})$, then $\tilde{u}_{k}$ satisfies
\begin{equation}
\partial_{y_i}(a_{ij}(\frac{y}{k})\partial_{y_j}\tilde{u}_{k}(y))-\frac{A}{|y|^{2}+1}\tilde{u}_{k}(y)=0,\quad y\in B(kx_{0},2kr)
\end{equation}
By Lemma \ref{Mean}, which applies to the elliptic case, we know that, for $y_{0}=kx_{0}$:
\begin{equation}
\tilde{u}_{k}^{2}(y_{0})\leq C\cdot\frac{\max\{\frac{kr}{1+k|x_{0}|},1\}^{-2\alpha}}{(kr)^{d}}\int_{B(y_{0},kr)}\tilde{u}_{k}^{2}(y)dy.
\end{equation}
Changing back to $x$-coordinate, we deduce
\begin{equation}
u_{k}^{2}(x_{0})\leq C\cdot\frac{\max\{\frac{r}{(1/k)+|x_{0}|},1\}^{-2\alpha}}{r^{d}}\int_{B(x_{0},r)}u_{k}^{2}(x)dx.
\end{equation}
By \eqref{limit}, since $r>|x_{0}|$ by assumption, after taking limit and using the dominated convergence theorem, we have:
\begin{equation}
u^{2}(x_{0})\leq C\cdot\Big(\frac{|x_{0}|}{r}\Big)^{2\alpha}\frac{1}{r^{d}}\int_{B(x_{0},r)}u^{2}(x)dx.
\end{equation}
Taking $r=1/2$, we know that $u$ decays to 0 as $x\to 0$ in an $\alpha-$H\"{o}lder sense.
\rightline{$\Box$}

We define $\Gamma_{a}(x,y)$ to be fundamental solution of the elliptic equation \eqref{maineq1} without the potential term, namely
\begin{equation}
\partial_{i}\big(a_{ij}(x)\partial_{j}\Gamma_{a}(x,y)\big)=\delta(x-y).
\end{equation}
Since the elliptic coefficients $a_{ij}$ are not smooth enough, we cannot hope for the gradient estimate of the related fundamental solution $\Gamma_{a}(x,y)$ as the Laplacian case. However, we have the following lemma on its H\"older estimate, which is a direct consequence of the De Giorgi-Nash-Moser theory. It plays a similar role as the gradient estimate in the proof of section 2.1, estimate of $I_{1}$ to $I_{3}$.

\begin{lemma}\label{levar}
$\Gamma_{a}(x,y)$ satisfies the following pointwise estimate
\begin{equation}\label{aest}
|\Gamma_{a}(x,y)|\leq C_{1} |x-y|^{2-d}
\end{equation}
and 
\begin{equation}\label{agest}
\begin{split}
|\Gamma_{a}(x_{1},y)-\Gamma_{a}(x_{2},y)|\leq C_{2}|x_{1}-x_{2}|^{\alpha}\Big[|x_{1}-y|^{2-d-\alpha}+&|x_{2}-y|^{2-d-\alpha} \Big]\\
&\forall y\in B(z,\delta)^{c}
\end{split}
\end{equation}
where $z=\frac{x_{1}+x_{2}}{2}$, $\delta=|x_{1}-x_{2}|$ as in section 2, $C_{1}$ and $C_{2}$ are two a priori constant depending only on $d$, $\lambda$, $\Lambda$.
\end{lemma}
\pf
The bound \eqref{aest} can be find in \cite{Aronson1967}, section 5. Thus we omit the proof here.\\
Now we prove \eqref{agest} for completeness. Consider $\Gamma_{a}(x,y)$ with $x\in B(z,\frac{3|y-z|}{4})$ and $y\in B(z,\delta)^{c}$. Since $\Gamma_{a}(x,y)$ is a well defined weak solution of the homogeneous elliptic equation
\begin{equation}\label{heq1}
\partial_{x_{i}}\big(a_{ij}(x)\partial_{x_{j}}\Gamma_{a}(x,y)\big)=0
\end{equation}
in domain $B(z,\frac{3|y-z|}{4})$, according to the classical H\"{o}lder estimate of weak solutions, (c.f.\cite{DeGiorgi1957}\cite{Moser1960}\cite{Nash1958}), we have:
\begin{equation}\label{est}
|\Gamma_{a}(x_{1},y)-\Gamma_{a}(x_{2},y)|\leq C \sup_{x\in B(z,\frac{3|y-z|}{4})}|\Gamma_{a}(x,y)|\cdot\Bigg(\frac{|x_{1}-x_{2}|}{d(x_{1},x_{2})}\Bigg)^{\alpha},
\end{equation}
here
\begin{equation}\label{d}
\begin{split}
d(x_{1},x_{2})&=\min\{dist(x_{1},\partial B(z,\frac{3|y-z|}{4})),dist(x_{2},\partial B(z,\frac{3|y-z|}{4}))\}\\
&=\frac{3|y-z|}{4}-\frac{\delta}{2}\geq C_{0}\max\{|y-x_{1}|,|y-x_{2}|\},\quad\forall y\in B(z,\delta)^{c}.
\end{split}
\end{equation}
Here $C_{0}$ is a constant. Now by \eqref{aest}
\begin{equation}\label{Gamma}
\begin{split}
&\sup_{x\in B(z,\frac{3|y-z|}{4})}|\Gamma_{a}(x,y)|\leq C_{1}\sup_{x\in B(z,\frac{3|y-z|}{4})}|x-y|^{2-d}\\
&\leq \tilde{C}\Big(|x_{1}-y|^{2-d}+|x_{2}-y|^{2-d}\Big)\quad\forall y\in B(z,\delta)^{c}
\end{split}
\end{equation}
where $\tilde{C}=\tilde{C}(d,\lambda,\Lambda)$.
By \eqref{est}, \eqref{d} and \eqref{Gamma}, we have $\forall y\in B(z,\delta)^{c}$
\begin{equation}
\begin{split}
&\quad|\Gamma_{a}(x_{1},y)-\Gamma_{a}(x_{2},y)|\\
&\leq\tilde{C}\Big(|x_{1}-y|^{2-d}+|x_{2}-y|^{2-d}\Big)\Big(C_{0}^{-1}|x_{1}-x_{2}|\Big)^{\alpha}\cdot\min\{|x_{1}-y|^{-\alpha},|x_{2}-y|^{-\alpha}\}\\
&\leq C_{2}|x_{1}-x_{2}|^{\alpha}\Big(|x_{1}-y|^{2-d-\alpha}+|x_{2}-y|^{2-d-\alpha}\Big)
\end{split}
\end{equation}
Thus we get \eqref{agest}.\\
\rightline{$\Box$}
Now we continue with the proof of the theorem. Suppose $u$ is a weak solution of \eqref{maineq1}, then it can be divided into
\begin{equation}\label{66}
u(x)=u_{a,0}(x)+\int_{B_{1/2}}\Gamma_{a}(x,y)\frac{A\cdot u(y)}{|y|^{2}}dy
\end{equation}
where $u_{a,0}(x)$ satisfies the homogeneous elliptic equation
\begin{equation}\label{heq}
\partial_{x_{i}}\big(a_{ij}(x)\partial_{x_{j}}u_{a,0}(x)\big)=0
\end{equation}
weakly in $B_{1/2}$. By the standard De Giorgi-Nash-Moser's theory, $u_{0}\in C^{\alpha_{0}}$, for some $\alpha_{0}=\alpha_{0}(d,\lambda,\Lambda)\in(0,1)$. Thus we only need to prove the H\"{o}lder continuity of
\begin{equation}
w_{a}(x)=\int_{B_{1/2}}\Gamma_{a}(x,y)\frac{A\cdot u(y)}{|y|^{2}}dy
\end{equation}
Similar to the Laplacian case in section 2.1, we have: using Lemma \ref{MeanValue} 
\begin{equation}
\begin{split}
|w_{a}(x_{1})-w_{a}(x_{2})|&\leq C\int_{B_{1/2}}|\Gamma_{a}(x_{1},y)-\Gamma_{a}(x_{2},y)|\frac{1}{|y|^{2-\alpha}}dy\\
&=C\int_{B_{1/2}\bigcap B(z,\delta)}|\Gamma_{a}(x_{1},y)|\frac{1}{|y|^{2-\alpha}}dy\\
&+C\int_{B_{1/2}\bigcap B(z,\delta)}|\Gamma_{a}(x_{2},y)|\frac{1}{|y|^{2-\alpha}}dy\\
&+C\int_{B_{1/2}- B(z,\delta)}|\Gamma_{a}(x_{1},y)-\Gamma_{a}(x_{2},y)|\frac{1}{|y|^{2-\alpha}}dy\\
&=C(I_{a,1}+I_{a,2}+I_{a,3})\\
\end{split}
\end{equation}
where $\alpha$ is identically same as in Lemma \ref{MeanValue}. By \eqref{aest} in Lemma \ref{levar}, we have the estimate of $I_{a,1}$ and $I_{a,2}$ is identically same as the estimate of $I_{1}$ and $I_{2}$ in section 3.1 since the fundamental solutions $\Gamma(x,y)$ and $\Gamma_{a}(x,y)$ share the same bound \eqref{aest}. Now we use estimate \eqref{agest} in Lemma \ref{levar} to estimate $I_{a,3}$, for $p<{d}/{(2-\alpha)}$:
\begin{equation}\label{eem}
\begin{split}
I_{a,3}&\leq C\delta^{\alpha}\int_{B_{1/2}- B(z,\delta)}\big(|x_{1}-y|^{2-\alpha-d}+|x_{2}-y|^{2-\alpha-d}\big)\frac{1}{|y|^{2-\alpha}}dy\\
& \leq C\delta^{\alpha}\Big(\int_{B(0,\frac{\delta}{2})^{c}}|y|^{\frac{(2-\alpha-d)p}{p-1}}dy\Big)^{1-1/p}\cdot\Big(\int_{B_{1/2}}|y|^{(\alpha-2)p}dy\Big)^{1/p}\\
&\leq C\delta^{2-\frac{d}{p}}\cdot\frac{(p-1)^{1-1/p}}{(\alpha-2)p+d}\\
&\leq C\cdot \delta^{\alpha_{-}},\quad \text{by choosing } p\to\big(\frac{d}{2-\alpha}\big)_{-}
\end{split}
\end{equation}
As mentioned before, we use $C_{-}$ to denote an arbitrary number close but smaller than $C$. Then we get the H\"{o}lder continuity of the potential $w_{a}(x)$, thus we finished the proof of the elliptic part of Theorem 2.
\section{Parabolic Case}
We consider \eqref{maineq2}
\begin{equation}\label{pbl}
\partial_{i} ( a_{ij}(x,t) \partial_{j}u(x,t)) - \frac{A}{|x|^{2}} u(x,t)-\partial_{t}u(x,t) =0
\end{equation}
in $B\times\mathbb{R}_{+}$, where $A>0$ is a constant and $a_{ij}(x,t)\in L^{\infty}$ satisfies the elliptic condition $\lambda I<\big(a_{ij}(x,t)\big)<\Lambda I$ with $0<\lambda\leq\Lambda<\infty$ as before. By classical De Giorgi's result, the solution is bounded on $B_{1/2}\times ]0,T]$ for some $0<T<\infty$. Similar to the elliptic case, treating the 0-order term $\frac{A}{|x|^{2}}u(t,x)$ as an external force term, and by Duhamel's Principle, for $t>t_{2}>0$, we could write a weak solution of \eqref{pbl} in the following form
\begin{equation}
\begin{split}
u(x,t)=u_{a,0}(x,t)+\int_{t_{2}/2}^{t}\int_{B_{1/2}}\Gamma_{a}(x,t;y,s)\frac{A\cdot u(y,s)}{|y|^{2}}dyds\\
\end{split}
\end{equation}
where $u_{a,0}$ is a weak solution of \eqref{pbl} in $B_{1/2}\times \mathbb{R}_{+}$ without potential ($A=0$) for $t>\frac{t_{2}}{2}$, namely
\begin{equation}\label{A=0}
\partial_{i} ( a_{ij}(x,t) \partial_{j}u(x,t))-\partial_{t}u(x,t) =0,
\end{equation}
and $\Gamma_{a}(x,t;y,s)$ is the fundamental solution of \eqref{A=0} with the source point at $(y,s)$. The first term is a weak solution of \eqref{A=0} and it is obviously H\"{o}lder continuous when $t>t_{2}/2>0$. Similarly as before, we only need to estimate
\begin{equation}
w_{a}(x,t)=\int_{t_{2}/2}^{t}\int_{B_{1/2}}\Gamma_{a}(x,t;y,s)\frac{A\cdot u(s,y)}{|y|^{2}}dyds
\end{equation}
and get its H\"{o}lder continuity in a space-time cube. Let $t_{1}>t_{2}>0$, $x_{1}\,,x_{2}\in B_{1/2}$, and set $z=\frac{x_{1}+x_{2}}{2}$, $\delta=|x_{1}-x_{2}|$ as before. We firstly give the parabolic version of Lemma \ref{MeanValue}
\begin{lemma}
Let u be a weak solution of the equation
\begin{equation}
\partial_{i}(a_{ij}(x,t)\partial_{j}u(x,t))-\frac{A}{|x|^{2}}u(x,t)-\partial_{t}u(x,t)=0
\end{equation}
$x_{0}\in B_{1/2}$, $\sqrt{t}/3\geq r\geq|x_{0}|\geq 0$. Here $a_{ij}\in L^{\infty}(B\times \mathbb{R}_{+})$, $\lambda I\leq \big(a_{ij}(x,t)\big)_{1 \leq i,j\leq d}\leq\Lambda I$ for some $0<\lambda\leq\Lambda<\infty$. Then there exists an $\alpha>0$, $C>0$ depending only on $\lambda$, $\Lambda$, $A$, $d$ and the $L^{2}$ module of $u$ on $Q_{r}(x_{0},t)$, such that
\begin{equation}\label{Pbound}
|u(x_{0},t)|^{2}\leq C t^{-(d/2+1+\alpha)}|x_{0}|^{2\alpha}
\end{equation}
\end{lemma}
 We could get this lemma similarly as the elliptic case (Lemma \ref{MeanValue}, the mean value inequality in section 3 is originally for parabolic case). We omit the details here.\\
\rightline{$\Box$}

We commence with the proof of the parabolic part of the Theorem 2.
\begin{equation}\label{MEstimate}
\begin{split}
&|w_{a}(x_{1},t_{1})-w_{a}(x_{2},t_{2})|\\
&=\bigg|\int_{\frac{t_{2}}{2}}^{t_{1}}\int_{B_{1/2}}\Gamma_{a}(x_{1},t_{1};y,s)\frac{A\cdot u(s,y)}{|y|^{2}}dyds\\
&\quad\quad-\int_{\frac{t_{2}}{2}}^{t_{2}}\int_{B_{1/2}}\Gamma_{a}(x_{2},t_{2};y,s)\frac{A\cdot u(s,y)}{|y|^{2}}dyds\bigg|\\
&\leq\int_{t_{2}}^{t_{1}}\int_{B_{1/2}}\Big|\Gamma_{a}(x_{1},y;t_{1},s)\Big|\cdot\Big|\frac{A\cdot u(s,y)}{|y|^{2}}\Big|dyds\\
&+\int_{\frac{t_{2}}{2}}^{t_{2}}\int_{B_{1/2}-B(z,\delta)}\Big|\Gamma_{a}(x_{1},t_{1};y,s)-\Gamma_{a}(x_{2},t_{2};y,s)\Big|\cdot\Big|\frac{A\cdot u(s,y)}{|y|^{2}}\Big|dyds\\
&+\int_{\frac{t_{2}}{2}}^{t_{2}}\int_{B(z,\delta)}\Big|\Gamma_{a}(x_{1},t_{1};y,s)\Big|\cdot\Big|\frac{A\cdot u(s,y)}{|y|^{2}}\Big|dyds\\
&+\int_{\frac{t_{2}}{2}}^{t_{2}}\int_{B(z,\delta)}\Big|\Gamma_{a}(x_{2},t_{2};y,s)\Big|\cdot\Big|\frac{A\cdot u(s,y)}{|y|^{2}}\Big|dyds\\
&:=I+II+III+IV.
\end{split}
\end{equation}


To bound term I, we choose $\alpha_{1}\in (0,\alpha)$, by the bound of $\Gamma_{a}$.We first bound the following intergal:
\begin{equation}
\begin{split}
&\int_{t_{2}}^{t_{1}}\frac{1}{(t_{1}-s)^{d/2}}\exp\Big(-\frac{|x_{1}-y|^{2}}{C(t_{1}-s)}\Big)s^{-(d/2+1+\alpha)}ds\\
&\leq t_{2}^{-(d/2+1+\alpha)}\int_{0}^{t_{1}-t_{2}}\frac{1}{s^{d/2}}\exp\Big(-\frac{|x_{1}-y|^{2}}{Cs}\Big)ds\\
&\leq t_{2}^{-(d/2+1+\alpha)}\cdot (t_{1}-t_{2})^{\alpha_{1}/2}\cdot\Big(\int_{0}^{t_{1}-t_{2}}\Big|\frac{1}{s^{d/2}}\exp\Big(\frac{-|x_{1}-y|^{2}}{Cs}\Big)\Big|^{\frac{2}{2-\alpha_{1}}}ds\Big)^{1-\alpha_{1}/2}\\
&\leq t_{2}^{-(d/2+1+\alpha)}\cdot(t_{1}-t_{2})^{\alpha_{1}/2}\cdot\Big(\int_{0}^{\infty}\Big|\frac{1}{s^{d/2}}\exp\Big(\frac{-|x_{1}-y|^{2}}{Cs}\Big)\Big|^{\frac{2}{2-\alpha_{1}}}ds\Big)^{1-\alpha_{1}/2}\\
&\lesssim t_{2}^{-(d/2+1+\alpha)}\cdot(t_{1}-t_{2})^{\alpha_{1}/2}\cdot|x_{1}-y|^{2-\alpha_{1}-d}.\\
\end{split}
\end{equation}
Thus, by \eqref{Pbound}, term I in \eqref{MEstimate} satisfies the following estimate:
\begin{equation}\label{E2}
I\lesssim t_{2}^{-(d/2+1+\alpha)}\cdot(t_{1}-t_{2})^{\alpha_{1}/2}\cdot\int_{B_{1/2}}|x_{1}-y|^{2-d-\alpha_{1}}\frac{1}{|y|^{2-\alpha}}dy.
\end{equation}
Choose $p>1$ such that $\frac{d}{2-\alpha_{1}}<p<\frac{d}{2-\alpha}$. By H\"{o}lder inequality, we have
\begin{equation}
\begin{split}
I&\lesssim t_{2}^{-(d/2+1+\alpha)}\cdot(t_{1}-t_{2})^{\alpha_{1}/2}\cdot \Big(\int_{B_{1/2}}|x_{1}-y|^{(2-d-\alpha_{1})\frac{p}{p-1}}dy\Big)^{1-1/p}\\
&\quad\cdot\Big(\int_{B_{1/2}}|y|^{p(\alpha-2)}dy\Big)^{1/p}\\
&\lesssim t_{2}^{-(d/2+1+\alpha)}\cdot(t_{1}-t_{2})^{\alpha_{1}/2}\cdot\Big(\int_{B}|y|^{(2-d-\alpha_{1})\frac{p}{p-1}}dy\Big)^{1-1/p}\\
&\quad\cdot\Big(\int_{B_{1/2}}|y|^{p(\alpha-2)}dy\Big)^{1/p}\\
&\lesssim t_{2}^{-(d/2+1+\alpha)}\cdot\Big(\frac{p-1}{(2-\alpha_{1})p-d}\Big)^{1-1/p}\cdot\Big(\frac{1}{(\alpha-2)p+d}\Big)^{1/p}\cdot(t_{1}-t_{2})^{\alpha_{1}/2}\\
&\lesssim t_{2}^{-(d/2+1+\alpha)}\cdot(t_{1}-t_{2})^{\frac{\alpha_{-}}{2}}.\hskip 1cm \text{as } \alpha_{1}\to\alpha_{-}
\end{split}
\end{equation}
Before the estimation of II, we first give two identity of the heat kernel function, which will be used later. Next we recall the following well-known facts for heat kernel.

Let
\begin{equation}\label{Ffunc}
f(t)=\frac{1}{t^{d/2}}\exp\Big(-\frac{a^{2}}{t}\Big),\quad a>0.
\end{equation}

Then
\begin{equation}\label{L1}
\int_{0}^{\infty}f(t)dt=\Gamma(d/2-1)a^{2-d}\sim a^{2-d},
\end{equation}
where $\Gamma$ is the usual $\Gamma$ function, and
\begin{equation}\label{Linfty}
\sup_{t\in(0,\infty)}f(t)=\exp(-\frac{d}{2})\cdot(\frac{d}{2})^{d/2}\cdot a^{-d}\sim a^{-d},
\end{equation}
and $f$ increases when $t\in(0,\frac{2a^{2}}{d})$, decreases when $t\in(\frac{2a^{2}}{d},\infty)$.

Since $y\in B(z,\delta)^{c}$, then for fixed $y$ and $s$, $\Gamma_{a}$ is a well-defined weak solution of equation
\begin{equation}\label{pblnp}
\partial_{x_{i}} ( a_{ij}(x,t) \partial_{x_{j}}u(x,t))-\partial_{t}u(x,t) =0
\end{equation}
for $(x,t)\in B(z,\frac{3|z-y|}{4})\times [t_{2}-\tau_{0},t_{1}]$, where $\tau_{0}$ is to be defined later. By the classical De Giorgi-Nash-Moser estimate, we have
\begin{equation}\label{two}
\begin{split}
&\Big|\Gamma_{a}(x_{1},t_{1};y,s)-\Gamma_{a}(x_{2},t_{2};y,s)\Big|\\
&\lesssim\sup_{(x,t)\in B(z,\frac{3|z-y|}{4})\times [t_{2}-\tau_{0},t_{1}]}\Big|\Gamma_{a}(x,t;y,s)\Big|\cdot\Big(\frac{|x_{1}-x_{2}|+\sqrt{t_{1}-t_{2}}}{\min\{|y-x_{1}|,|y-x_{2}|,\sqrt{\tau_{0}}\}}\Big)^{\alpha}\\
&\lesssim \sup_{(x,t)\in B(z,\frac{3|z-y|}{4})\times [t_{2}-\tau_{0},t_{1}]}\frac{1}{(t-s)^{d/2}}\exp(-\frac{|x-y|^{2}}{C(t-s)})\cdot\Big(\frac{|x_{1}-x_{2}|+\sqrt{t_{1}-t_{2}}}{\min\{|y-x_{1}|,|y-x_{2}|,\sqrt{\tau_{0}}\}}\Big)^{\alpha}\\
&\lesssim\sup_{t\in[t_{2}-\tau_{0},t_{1}]}\frac{1}{(t-s)^{d/2}}\Big(\exp(-\frac{|x_{1}-y|^{2}}{C(t-s)})+\exp(-\frac{|x_{2}-y|^{2}}{C(t-s)})\Big)\cdot\Big(\frac{|x_{1}-x_{2}|+\sqrt{t_{1}-t_{2}}}{\min\{|y-x_{1}|,|y-x_{2}|,\sqrt{\tau_{0}}\}}\Big)^{\alpha}\\
\end{split}
\end{equation}
The second "$\lesssim$" is due to the bound of fundamental solution of parabolic equation, which can be found in \cite{Aronson1967}. Now choose $\tau_{0}=\min\{|x_{1}-y|^{2},|x_{2}-y|^{2}\}/C$, where $C$ is the constant in the last line of \eqref{two}. Since in the domain $y\in B(z,\delta)^{c}$, $\Gamma_{a}(x,t;y,s)$ could do a well-defined 0 extension to $t<s$, we need not worry about $t_{2}<\tau_{0}$. In order to proceed, we need to bound the following integral:
\begin{equation}
\begin{split}
J:=&\int_{\frac{t_{2}}{2}}^{t_{2}}\sup_{t\in[t_{2}-\tau_{0},t_{1}]}\frac{1}{(t-s)^{d/2}}\Big(\exp(-\frac{|x_{1}-y|^{2}}{C(t-s)})+\exp(-\frac{|x_{2}-y|^{2}}{C(t-s)})\Big)\cdot s^{-(d/2+1+\alpha)}ds\\
\lesssim& t_{2}^{-(d/2+1+\alpha)}\Bigg[\int_{\frac{t_{2}}{2}}^{t_{2}}\sup_{t\in[t_{2}-\tau_{0},t_{1}]}\frac{1}{(t-s)^{d/2}}\Big(\exp(-\frac{|x_{1}-y|^{2}}{C(t-s)}))\Big)ds\\
&+\int_{\frac{t_{2}}{2}}^{t_{2}}\sup_{t\in[t_{2}-\tau_{0},t_{1}]}\frac{1}{(t-s)^{d/2}}\Big(\exp(-\frac{|x_{2}-y|^{2}}{C(t-s)})\Big)ds\Bigg]\\
:=&t_{2}^{-(d/2+1+\alpha)}\cdot (J_{1}+J_{2}).
\end{split}
\end{equation}
Note that this integral has a "sup" inside, we will need to split the interval of integration here. Let $C_{0}=\frac{1}{C}\big(1+\frac{2}{d}\big)$, for $i=1,2$, we have, when $s\in [t_{2}/2,t_{2}-C_{0}|x_{i}-y|^{2}]$ and $t\in [t_{2}-\tau_{0},t_{1}]$, $t-s\geq\frac{2}{d}\cdot\frac{|x_{i}-y|^{2}}{C}$.

If $t_{2}>2C_{0}|x_{i}-y|^{2}$, namely $t_{2}-C_{0}|x_{i}-y|^{2}>t_{2}/2$, we have
\begin{equation}
\begin{split}
J_{i}\leq& \int_{\frac{t_{2}}{2}}^{t_{2}-C_{0}|x_{i}-y|^{2}}\sup_{t\in[t_{2}-\tau_{0},t_{1}]}\frac{1}{(t-s)^{d/2}}\exp(-\frac{|x_{i}-y|^{2}}{C(t-s)})ds\\
+&\int_{t_{2}-C_{0}|x_{i}-y|^{2}}^{t_{2}}\sup_{t\in[t_{2}-\tau_{0},t_{1}]}\frac{1}{(t-s)^{d/2}}\exp(-\frac{|x_{i}-y|^{2}}{C(t-s)})ds\\
\end{split}
\end{equation}
By \eqref{Ffunc}, the monotonicity of $f$, we have
\begin{equation}
\begin{split}
\sup_{t\in[t_{2}-\tau_{0},t_{1}]}\frac{1}{(t-s)^{d/2}}\exp(-\frac{|x_{i}-y|^{2}}{C(t-s)})=\frac{1}{(t_{2}-\tau_{0}-s)^{d/2}}\exp(-\frac{|x_{i}-y|^{2}}{C(t_{2}-\tau_{0}-s)})\\
\forall s\in [t_{2}/2,t_{2}-C_{0}|x_{2}-y|^{2}]\hskip 0.5cm i=1,2.
\end{split}
\end{equation}

Thus, by \eqref{L1} and \eqref{Linfty}, we have:
\begin{equation}
\begin{split}
J_{i}&\lesssim \int_{0}^{\infty}\frac{1}{s^{d/2}}\exp(-\frac{|x_{i}-y|^{2}}{Cs})ds+C_{0}|x_{i}-y|^{2-d}\\
&\lesssim|x_{i}-y|^{2-d}.\\
\end{split}
\end{equation}
Else, if $t_{2}\leq 2C_{0}|x_{i}-y|^{2}$, we have:
\begin{equation}
\begin{split}
J_{i}\lesssim & |x_{i}-y|^{2}\cdot\sup_{s\in\mathbb{R}_{+}}\frac{1}{s^{d/2}}\exp\Big(-\frac{|x_{i}-y|^{2}}{Cs}\Big)\\
\lesssim & |x_{i}-y|^{2-d}.
\end{split}
\end{equation}
Therefore, term II satisfies the following estimate
\begin{equation}
\begin{split}
II&\lesssim t_{2}^{-(d/2+1+\alpha)}\int_{B_{1/2}-B(z,\delta)}\frac{J_1+J_2}{|y|^{2-\alpha}}dy\\
&\lesssim t_{2}^{-(d/2+1+\alpha)}\cdot(|x_{1}-x_{2}|+\sqrt{t_{1}-t_{2}})^{\alpha}\\
&\hskip 1cm\cdot\int_{B_{1/2}-B(z,\delta)}\big(|x_{1}-y|^{2-\alpha-d}+|x_{2}-y|^{2-\alpha-d}\big)\frac{1}{|y|^{2-\alpha}}dy\\
&\lesssim t_{2}^{-(d/2+1+\alpha)}\cdot (|x_{1}-x_{2}|+\sqrt{t_{1}-t_{2}})^{\alpha}\\
&\hskip 1cm\cdot\Big(\int_{B(0,\frac{\delta}{2})^{c}}|y|^{\frac{(2-\alpha-d)p}{p-1}}dy\Big)^{1-1/p}\cdot\Big(\int_{B_{1/2}}|y|^{(\alpha-2)p}dy\Big)^{1/p}\\
&\lesssim t_{2}^{-(d/2+1+\alpha)}\cdot(|x_{1}-x_{2}|+\sqrt{t_{1}-t_{2}})^{2-\frac{d}{p}}\cdot\frac{(p-1)^{1-1/p}}{(\alpha-2)p+d}\\
&\lesssim t_{2}^{-(d/2+1+\alpha)}\cdot(|x_{1}-x_{2}|+\sqrt{t_{1}-t_{2}})^{\alpha_{-}},\quad \text{by choosing}\,\,p\to\big(\frac{d}{2-\alpha}\big)_{-}.
\end{split}
\end{equation}
Finally, III and IV are essentially the same, thus we estimate them together
\begin{equation}\label{(5.22)}
\begin{split}
III,IV&\lesssim t_{2}^{-(d/2+1+\alpha)}\cdot\int_{0}^{\infty}\int_{B(z,\delta)}\frac{1}{s^{d/2}}\exp(-\frac{|x_{i}-y|^{2}}{Cs})\Big|\frac{1}{|y|^{2-\alpha}}\Big|dyds,\quad i=1,2\\
&\lesssim t_{2}^{-(d/2+1+\alpha)}\cdot\int_{B(z,\delta)}|x_{i}-y|^{2-d}\frac{1}{|y|^{2-\alpha}}dy,\quad i=1,2\\
&\lesssim t_{2}^{-(d/2+1+\alpha)}\cdot|x_{1}-x_{2}|^{\alpha_{-}}.\\
\end{split}
\end{equation}
The last step is the same as \eqref{I1} in section 3. Thus we get the H\"{o}lder continuity of $w_{a}(t,x)$ with $x_{1},\,x_{2}\in B_{1/2}$ and $t_{1}>t_{2}>0$. The H\"{o}lder norm of $u$ respect to $\lambda$, $\Lambda$, $d$, $A$, $t_{2}$, and the $L^{2}$ norm of $u$, i.e.
\begin{equation}
\|u\|_{C^{\alpha_{-};\frac{\alpha_{-}}{2}}([t_{0},\infty)\times B_{1/2})}\leq C(\lambda,\Lambda,d,A)t_{0}^{-(d/2+1+\alpha)}\|u\|_{L^{2}}
\end{equation}
where $\alpha=\alpha(\lambda,\Lambda,d,A)>0$.
\section{Appendix}

\subsection{Existence Results}
Since the Potencial $\frac{A}{|x|^{2+\beta}}$ is much more singular than the classical case, the existence and uniqueness of weak solutions of \eqref{maineq2} can not be taken for granted. Let $u_{k}$ be the unique solution to the following parabolic problem:
\begin{equation}\label{APP}
\begin{split}
\pa_{i}\big(a_{ij}\pa_{j}u_{k}\big)-Ac_{k}u_{k}-\pa_{t}u_{k}&=0.\hskip 0.5cm \text{in } B\times [0,T]\\
u_{k}(x,0)&=u_{0}\in L^{2}(B)\hskip 0.5cm \text{on  }B\\
u_{k}&=0\hskip 0.5cm \text{on  }\pa B\times[0,T].
\end{split}
\end{equation}
Here $a_{ij}$ is as in \eqref{maineq2}, $c_{k}=c_{k}(x)\in C^{\infty}_{0}(B)$ satisfies
\begin{equation}
c_{k}\to\frac{1}{|x|^{2+\beta}},\quad\text{strongly in }L^{p}(B).
\end{equation}
for a fixed $p\in[1,\frac{d}{2+\beta}[$. This is possible since $\frac{1}{|x|^{2+\beta}}\in L^{(\frac{d}{2+\beta})_{-}}(B)$. Multiplying Eq.\eqref{APP} by $u_{k}$, one easily obtain
\begin{equation}
\lambda\int_{0}^{T}\int_{B}|\nabla u_{k}|^{2}dxdt+A\int_{0}^{T}\int_{B}c_{k}u_k^{2}dxdt+\int_{B}u_k^{2}(\cdot,T)dx\leq\int_{B}u^{2}_{0}dx.
\end{equation}
Hence there exists a function $u$ such that $u$, $\nabla u\in L^{2}(B\times[0,T])$ and a subsequence of $\{u_k\}$ such that
\begin{equation}
\begin{split}
u_{k}&\rightharpoonup u\quad\text{weakly in }L^{2}(B\times[0,T]);\\
\nabla u_{k}&\rightharpoonup \nabla u\quad\text{weakly in }L^{2}(B\times[0,T]);\\
u_{k}&\rightharpoonup u\quad\text{weakly in }L^{2}([0,T],L^{\frac{2d}{d-2}}(B)).
\end{split}
\end{equation}
Now we are going to prove $u$ is a weak solution of \eqref{maineq2} with the same initial and boundary condition as \eqref{APP}. Clearly, $\forall\phi\in C^{\infty}_{0}(B\times[0,T])$, $u_{k}$ satisfies
\begin{equation}\label{APPW}
\int_{0}^{T}\int_{B}a_{ij}\pa_{i}u_{k}\pa_{j}\phi dxdt+A\int_{0}^{T}\int_{B}c_{k}u_{k}\phi dxdt-\int_{0}^{T}\int_{B}u_{k}\pa_{t}\phi dxdt=\int_{B}u_{0}\phi(x,0) dx.
\end{equation}
By the weak convergence of $u_{k}$ and $\nabla u_{k}$, we have
\begin{equation}
\begin{split}
\int_{0}^{T}\int_{B}a_{ij}\pa_{i}u_{k}\pa_{j}\phi dxdt-\int_{0}^{T}\int_{B}u_{k}\pa_{t}\phi dxdt\to&\int_{0}^{T}\int_{B}a_{ij}\pa_{i}u\pa_{j}\phi dxdt-\int_{0}^{T}\int_{B}u\pa_{t}\phi dxdt,\\
&\hskip 3cm\text{as }k\to\infty.
\end{split}
\end{equation}
Next, notice that
\begin{equation}
\begin{split}
&\int_{0}^{T}\int_{B}c_{k}u_{k}\phi dxdt-\int_{0}^{T}\int_{B}\frac{u\phi}{|x|^{2+\beta}}dxdt\\
=&\int_{0}^{T}\int_{B}u_{k}\phi\Big(c_{k}-\frac{1}{|x|^{2+\beta}}\Big)dxdt+\int_{0}^{T}\int_{B}\frac{1}{|x|^{2+\beta}}(u_{k}-u)\phi dxdt.\end{split}
\end{equation}
For $\beta\in[0,d-2[$, by the strong convergence of $c_{k}$ and weak convergence of $u_{k}$, we have
\begin{equation}\label{LC}
\int_{0}^{T}\int_{B}c_{k}u_{k}\phi dxdt-\int_{0}^{T}\int_{B}\frac{u\phi}{|x|^{2+\beta}}dxdt\to 0\quad\text{as }k\to\infty.
\end{equation}
By \eqref{APPW} and \eqref{LC}, we obtain
\begin{equation}
\int_{0}^{T}\int_{B}a_{ij}\pa_{i}u\pa_{j}\phi dxdt+A\int_{0}^{T}\int_{B}\frac{1}{|x|^{2+\beta}}u\phi dxdt-\int_{0}^{T}\int_{B}u\pa_{t}\phi dxdt=\int_{B}u_{0}\phi(x,0) dx.
\end{equation}
i.e. $u$ is a weak solution to \eqref{maineq2}. Since elliptic problem is a time-independent parabolic case, we get the existence results for both.

\subsection{Local Boundedness of the Weak Solution and Maximum Principle}
\begin{lemma}\label{Bdd}
The weak solution of \eqref{maineq1} in $B$ is bounded in $B_{1/2}$.
\end{lemma}
\pf
Since the potential term coefficient $\frac{A}{|x|^{2+\beta}}\geq 0$, we could get boundedness of $u$ by classical De Giorgi iteration method. We omit the details here.
\rightline{$\Box$}

\begin{lemma}[weak maximum principle for elliptic equation]\label{Wmp}
Let $u\in\mathcal{H}^{\beta}(B)$ be a weak solution of elliptic equation \eqref{maineq1} in $B$, then
\begin{equation}
\sup_{B_{1/2}}u\leq\sup_{\partial B_{1/2}}u_{+},\quad\inf_{B_{1/2}}u\geq\inf_{\partial B_{1/2}}u_{-}
\end{equation}
Here $u_+=\max\{u,0\}$ and $u_-=\min\{u,0\}$.
\end{lemma}
\pf
Choose test function $v=max\{u-l,0\}$, where $l=\sup_{B_{1/2}}u_{+}$, pay attention that $v\in \mathcal{H}^{\beta}_{0}(B_{1/2})$, $\frac{1}{|x|^{2+\beta}}u\cdot v\geq 0$, then we have
\begin{equation}
\begin{split}
0&\geq\int_{B_{1/2}}a_{ij}(x)\partial_{i}u(x)\partial_{j}v(x)dx\\
&=\int_{B_{1/2}}a_{ij}(x)\partial_{i}(u(x)-l)_{+}\partial_{j}(u(x)-l)_{+}dx\\
&\geq\lambda\int_{B_{1/2}}|\nabla(u-l)_{+}|^{2}dx
\end{split}
\end{equation}
which means $v\equiv 0$ in $B_{1/2}$, which we get the first inequality. The second one is similar.\\
\rightline{$\Box$}
\subsection{An Introduction to the Modified Bassel's Equation}
We call an important ordinary differential equation which was used in this paper: \emph{the Modified Bessel's equation}
\begin{equation}\label{MBessel}
t^{2}x''(t)+tx'(t)-(t^{2}+\lambda^{2})x(t)=0.
\end{equation}
This equation has two linearly independent solution,i.e.
$\mathcal{I}_{\lambda}(t)$, $\mathcal{K}_{\lambda}(t)$ ,
which are exponentially growing and decaying as $t\to+\infty$ and which are referred to
as \emph{modified Bessel's function}of first and second kind, respectively.
For more detailed information of \emph{Modified Bessel's equation} and
\emph{modified Bessel's function}, we refer to \cite{Abramowitz1972},
a well-known handbook of mathematical special functions.\\
\medskip

{\bf Acknowledgment.} Q. S. Z. gratefully acknowledges
the support by Simons foundation and Siyuan Foundation through Nanjing University.


\section{References}
\begin{thebibliography}{10}

\bibitem{Abramowitz1972}
Milton Abramowitz and Irene~A. Stegun.
\newblock {\em Handbook of mathematical functions with formulas, graphs, and
  mathematical tables}.
\newblock Dover Publications, Inc., 1992 edition, 1992.

\bibitem{Aronson1967}
D.~G. Aronson.
\newblock Bounds for the fundamental solution of a parabolic equation.
\newblock {\em Bulletin of the American Mathematical Society}, 73:890--896,
  1967.

\bibitem{Baras1984}
Pierre Baras and Jerome~A. Goldstein.
\newblock The heat equation with a singular potential.
\newblock {\em Trans. Amer. Math. Soc.}, 284(1):121--139, 1984.

\bibitem{GT1998}
David Gilbarg and Neil~S. Trudinger.
\newblock {\em Elliptic Partial Differential Equations of Second Order}.
\newblock Springer, 2nd edition, 1998.

\bibitem{DeGiorgi1957}
E.~De Giorgi.
\newblock Sulla differenziabilit\`{a} e l'analiticit\`{a} delle estremali degli
  intergrali multipli regolari.
\newblock {\em Men. Accad. Sci. Torino. Cl. Sci. Fis. Mat. Natur}, (3):25--43,
  3 1957.

\bibitem{Alexander2006}
Alexander Grigor'yan.
\newblock {\em Heat kernels on weighted manifolds and applications.}, volume
  398 of {\em Contemp. Math.}
\newblock Amer. Math. Soc., 2006.

\bibitem{Milman2003}
Pierre~D. Milman and Yu.~A. Semenov.
\newblock Heat kernel bounds and desingularizing weights.
\newblock {\em J. Funct. Anal.}, 202(1):1--24, 2003.

\bibitem{Milman2004}
Pierre~D. Milman and Yu.~A. Semenov.
\newblock Global heat kernel bounds via desingularizing weights.
\newblock {\em J. Funct. Anal.}, 212(2):373--398, 2004.

\bibitem{Moser1960}
J.~Moser.
\newblock A new proof of de Giorgi's theorem concerning the regularity problem
  for elliptic differential equations.
\newblock {\em Comm. Pure Appl. Math}, 13:457--468, 1960.

\bibitem{Moser1961}
J.~Moser.
\newblock On Harnack's theorem for elliptic differential equations.
\newblock {\em Comm. Pure Appl. Math}, 14:577--591, 1961.

\bibitem{Nash1958}
J.~Nash.
\newblock Continuity of solutions of parabolic and elliptic equations.
\newblock {\em Amer. J. Math}, (80):931--954, 1958.

\bibitem{Wong2003}
B.~Wong and Qi~S. Zhang.
\newblock Refined gradient bounds, possion equations and some applications to
  open k\"{a}hler manifolds.
\newblock {\em Asian J. Math}, 7(3):1--28, September 2003.

\bibitem{Zhang2000}
Qi~S. Zhang.
\newblock Large time behavior of Schr\"{o}dinger heat kernels and applications.
\newblock {\em Comm. Math. Phys.}, 210(2):371--398, 2000.

\end{thebibliography}
\end{document}